\documentclass[12pt]{article}

\usepackage{amsthm}
\usepackage{mathtools}
\usepackage{hslTR}
\usepackage{rgsEnvironments}
\usepackage{rgsMacros}
\usepackage{amsmath,amssymb,amsfonts}
\usepackage{cite}
\usepackage{enumitem,kantlipsum}
\usepackage{hyperref}
\usepackage{etoolbox}
\usepackage{xfrac}
\usepackage[numbers,sort]{natbib}
\providetoggle{long}
\settoggle{long}{true}

\providetoggle{dissertation}
\settoggle{dissertation}{false}

\providetoggle{report}
\settoggle{report}{true}

\providetoggle{test}
\settoggle{test}{false}

\newcommand{\dobib}{ 
    \bibliographystyle{unsrt}
    \bibliography{Bibliography} 
}

\def\startmodifnew{%
\begingroup\color{blue}
}

\def\stopmodifnew{%
\protect\normalfont
\endgroup%
}

\newcommand*\squeezespaces[1]{
  \thickmuskip=\scalemuskip{\thickmuskip}{#1}%
  \medmuskip=\scalemuskip{\medmuskip}{#1}%
  \thinmuskip=\scalemuskip{\thinmuskip}{#1}%
  \nulldelimiterspace=#1\nulldelimiterspace
  \scriptspace=#1\scriptspace
}
\newcommand*\scalemuskip[2]{%
  \muexpr #1*\numexpr\dimexpr#2pt\relax\relax/65536\relax
}

\makeatletter
\newcommand*{\coleq}{\mathrel{\rlap{%
                     \raisebox{0.3ex}{$\m@th\cdot$}}%
                     \raisebox{-0.3ex}{$\m@th\cdot$}}%
                     =}
\makeatother

\AtBeginDocument{%
  \mathchardef\stdcomma=\mathcode`,
  \mathcode`,="8000
}
\begingroup\lccode`~=`, \lowercase{\endgroup\def~}{\stdcomma\,}

\begin{document}
\renewcommand{\dobib}{} 





\title{Robust Finite-time Parameter Estimation for Nonlinear Systems with Linear Parametric Uncertainty}

\author{Ryan S. Johnson, Adnane Saoud, and Ricardo G. Sanfelice}
\maketitle

\begin{abstract}
%
We consider the problem of estimating a vector of unknown constant parameters for a class of continuous-time systems with nonlinear dynamics that depend linearly on the unknown parameters. 
We propose a hybrid estimation algorithm -- that is, an algorithm whose state variables exhibit both continuous (flow) and discrete (jump) evolution -- that, under excitation over only a finite interval, guarantees finite-time convergence of the parameter estimate to the true value. 
Furthermore, we show that the parameter estimate is robust with respect to noise in the measurements of the plant state. 
Simulation results including a spacecraft application show the merits of the proposed approach.
\vspace{-7pt} \noindent
\end{abstract}

\section{Introduction}

\subsection{Background}
Estimating the unknown parameters of a system is crucial in many engineering applications \cite{8855027,8466598,doi:10.2514/3.20454}. 
Parameter estimation algorithms typically require system inputs to satisfy a persistence of excitation (PE) condition to ensure convergence of the parameter estimate to the true value \cite{Narendra_adaptive_1989,Sastry1989Adaptive,Khalil2002Nonlinear}. 
Namely, a continuous-time bounded signal $t \mapsto \phi(t) \in \reals^{n \times p}$ is persistently exciting if there exist $T, \mu_\phi > 0$ such that\footnote{The classical PE condition is defined for vector signals but it extends naturally to matrix signals. Since this paper considers a matrix regressor, we restate the PE condition using a matrix signal for ease of comparison.}
\begin{equation} \label{eqn:FTcontinuousPE}
    \int_{t'}^{t'+T} 
    \phi (t)^\top 
    \phi (t) dt
    \geq \mu_\phi I
    \quad \forall t' \geq 0.
\end{equation}
This PE condition is often difficult to verify online since the property must hold for all future time. 
Much work has been done to develop estimation algorithms that ensure convergence of the parameter estimate without requiring PE. 
Finite-time parameter estimators are one such class of algorithms.  
Since these algorithms achieve finite-time convergence of the parameter estimate to the true value, they may relax the PE condition to instead require excitation over only a finite interval -- henceforth called finite excitation (FE). 
Namely, given $t' \geq 0$ and $T > 0$, a continuous-time bounded signal $t \mapsto \phi(t) \in \reals^{n \times p}$ is exciting over the interval $[t', t' + T]$ if there exist $\mu_\phi > 0$ such that
\begin{equation} \label{eqn:FTcontinuousFE}
    \int_{t'}^{t'+T} 
    \phi (t)^\top 
    \phi (t) dt
    \geq \mu_\phi I.
\end{equation}
%
%
Finite-time parameter estimation for linear regression models was studied in \cite{243}, which employs a hybrid systems framework \cite{65, 220}, in \cite{Wang2019FiniteTime}, which employs dynamic regressor extension and mixing (DREM), and in \cite{Ortega1988LeastSquares}, which utilizes a singular value decomposition of the inverse covariance matrix. All of these works ensure estimation of unknown parameters under FE conditions similar to \eqref{eqn:FTcontinuousFE}.
Finite-time parameter estimation for classes of nonlinear systems was studied in \cite{Adetola2007FiniteTime}, which employs coupled state and parameter estimators, and in \cite{Na2015FiniteTime}, which employs sliding mode techniques. Both of these works impose a PE condition similar to \eqref{eqn:FTcontinuousPE}.
The work \cite{Vahidi2021Memory} proposes a memory-augmented finite-time estimation algorithm that ensures convergence of the estimate under a FE condition similar to \eqref{eqn:FTcontinuousFE} by leveraging stored samples of the system trajectory.
Several DREM-based approaches have recently been proposed. The work \cite{Mohamed2026DREM} proposes a finite-time estimation scheme for a class of nonlinear systems, but provides no guarantees on the numerical conditioning of the parameter update. The work \cite{Ortega2022ParamId} addresses the alertness problem (the loss of parameter tracking capability after convergence of the estimate) through a sliding window approach that maintains tracking capability.

The hybrid systems framework in \cite{65, 220} has also been used for parameter estimation in \cite{Saoud2024HybridPersistency}, which extends \cite{242}. Although \cite{Saoud2024HybridPersistency} and this paper rely on the same formalism and both address a parameter estimation problem, they differ in two fundamental ways: in \cite{Saoud2024HybridPersistency}, the signals involved in the estimation problem are hybrid and the parameter estimate converges asymptotically (in fact, exponentially), whereas here the plant dynamics are continuous and the parameter estimate converges in finite time. More precisely, in \cite{Saoud2024HybridPersistency} the hybrid behavior originates from the data: the regressor, or the plant generating it, exhibits both flows and jumps, and the estimator is a hybrid gradient algorithm whose jumps are synchronized with those of the data. The main contribution of \cite{Saoud2024HybridPersistency} is a hybrid notion of persistence of excitation which, when it holds over every window of the hybrid time domain, implies uniform exponential stability and input-to-state stability of the estimation error dynamics. In contrast, the plant considered here evolves in continuous time, and the hybrid behavior of the estimator is a design choice: the estimator jumps opportunistically, only when the excitation is sufficient to reconstruct the unknown parameter vector, which allows the estimate to be reset to the true value in finite time. The conditions and guarantees differ accordingly: the excitation condition in \cite{Saoud2024HybridPersistency} must hold for all hybrid time, whereas we require excitation over only a finite interval.

To the best of our knowledge, a hybrid algorithm that guarantees finite-time estimation of unknown parameters under excitation over only a finite interval has not been established in the literature for the considered class of systems. The goal of this paper is to fill that gap.

\subsection{Contributions}
In this paper, we propose a hybrid algorithm for finite-time estimation of unknown constant parameters in a class of continuous-time systems with nonlinear dynamics that are affine in the unknown parameters. 
We establish sufficient conditions that guarantee finite-time convergence of the parameter estimate to the true value, and we upper bound the convergence time of the parameter estimate. 
The main contributions in this paper are the following:
\begin{itemize}
    \item \underline{Estimation under finite excitation:} in the main stability result (Theorem~\ref{thm:FTClosed_Loop_Stab_HS1}), we establish that our algorithm guarantees finite-time convergence of the parameter estimate to the true value when the regressor is exciting over only a finite interval as in \eqref{eqn:FTcontinuousFE} and is Lipschitz continuous.
    Moreover, we ensure that the parameter estimate is numerically well conditioned at the convergence time by upper bounding the condition number of a matrix whose inverse is used to reset the parameter estimate at jumps.
    To the best of our knowledge, this is the first hybrid algorithm that ensures finite-time estimation of unknown parameters for the considered class of systems under excitation over only a finite interval. Furthermore, it is the first such algorithm that ensures the parameter estimate is numerically well conditioned at the convergence time.
    \item \underline{Robustness to state measurement noise:} we establish two main robustness results. The first result (Theorem~\ref{thm:FTISS}) establishes that, when the regressor is persistently exciting as in \eqref{eqn:FTcontinuousPE} and is Lipschitz continuous, our algorithm is input-to-state stable with respect to noise in the measurements of the plant state. The second result (Theorem~\ref{thm:FTiISS}) establishes that, under excitation over a finite interval, the parameter estimate is bounded by a function of the integral of the measurement noise.
\end{itemize}

In \cite{256}, we proposed a hybrid algorithm for estimating unknown constant parameters in a class of continuous-time systems with linear dynamics. 
%
%
This paper differs from \cite{256} in several key ways. 
Namely, here we consider a wider class of continuous-time systems with nonlinear dynamics and we propose a new hybrid algorithm to solve the estimation problem. We show that the proposed algorithm guarantees estimation of the unknown parameters under excitation over only a finite interval, which relaxes the condition imposed in \cite{256}. The proposed algorithm guarantees that the jump map is numerically well conditioned at jumps, whereas \cite{256} provides no such guarantee (see Remark~\ref{rem:FTConditionNumber}). This paper provides new robustness results for the proposed algorithm under both PE and FE conditions. Finally, new examples including a spacecraft application are presented.

Our approach is related to the ones in \cite{243,SaoudFiniteTimeInPrep,hartman2012robust} -- each proposes a finite-time estimator using a hybrid system framework. In comparison to our approach, 
the works \cite{243, SaoudFiniteTimeInPrep} consider a linear regression model while we consider a dynamical system. 
Moreover, the estimation algorithm in \cite{243, SaoudFiniteTimeInPrep} uses two coupled parameter estimates, while our algorithm uses only a single parameter estimate along with a memory state, which allows us to relax the conditions required for finite-time convergence of the estimate.
Finally, the algorithm in \cite{243, SaoudFiniteTimeInPrep} relies on a timer to trigger jumps in the estimator state, whereas our algorithm jumps opportunistically when reconstruction of the unknown parameter vector is possible.
In \cite{hartman2012robust}, the authors rely on a PE condition similar to \eqref{eqn:FTcontinuousPE} to ensure that their algorithm is well defined and to guarantee completeness of maximal solutions, 
while our algorithm requires excitation over only a finite interval. Moreover, the algorithm in \cite{hartman2012robust} jumps based on the determinant of a matrix, which does not guarantee that the matrix is numerically well conditioned for inversion at jumps (see Remark~\ref{rem:FTConditionNumber}). In contrast, our algorithm guarantees that the jump map is numerically well conditioned at jumps.

\subsection{Organization of the Paper}
The remainder of this paper is organized as follows. Preliminaries on hybrid systems are presented in Section~\ref{sec:FTPreliminaries}. Our hybrid algorithm is presented in Section~\ref{sec:FTUniningControl}, and the convergence and robustness properties are analyzed in Sections~\ref{sec:FTConvergenceAnalysis} and~\ref{sec:FTRobustnessAnalysis}, respectively. Simulation results are in Sections~\ref{sec:FTConvergenceAnalysis},~\ref{sec:FTRobustnessAnalysis}, and~\ref{sec:FTNumericalExample}. Conclusions and future work are in Section~\ref{sec:FTConclusion}. 
\iftoggle{report}{
The proofs of Lemmas~\ref{lem:FT_FEPsi},~\ref{lem:FTHSnuSol}, and~\ref{lem:FTepsISS} are in the Appendix.
}{
The proof of Lemma~\ref{lem:FT_FEPsi} is in the Appendix. Due to space constraints, the proofs of Lemmas~\ref{lem:FTHSnuSol} and~\ref{lem:FTepsISS} are in \cite{JohnsonFiniteTimeArxiv}.
}

\section{Preliminaries} \label{sec:FTPreliminaries}
%
\subsection{Notation}

We denote the set of real, nonnegative, positive, and natural numbers (including zero) as $\reals$, $\realsgeq$, $\realsg$, and $\nats$, respectively. The matrix $I$ denotes the identity matrix of appropriate dimension. 
%
%
For $x,y \in \reals^n$, we write $[x^\top \ \, y^\top]^\top$ as $(x,y)$.
The Euclidean norm of vectors and the associated induced matrix norm is denoted by $|\cdot|$.
The Frobenius norm of a matrix $M \in \reals^{n \times p}$ is $|M|_{\mathrm F} \coleq \sqrt{\tr (M^\top M)}$, where $\tr(\cdot)$ is the trace.
Given a matrix $A \in \reals^{n \times n}$, $\text{\rm eig} (A)$ denotes the set of all eigenvalues of $A$, $\lambda_{\max}(A) \coleq \max \{ \lambda/2 : \lambda \in \text{\rm eig} (A + A^\top)\}$ and $\lambda_{\min}(A) \coleq \min \{ \lambda/2 : \lambda \in \text{\rm eig} (A + A^\top)\}$. 
%
%
The condition number of a nonsingular matrix $A \in \reals^{n \times n}$ is $\text{\rm cond}(A) \coleq |A||A^{-1}|$.
The distance of a point $x$ to a nonempty set $S$ is denoted by $|x|_S = \inf_{y \in S} |y-x|$ (the quantity $|\cdot|_S$ should not be confused with the Frobenius norm $|\cdot|_{\mathrm F}$ since, besides the subscript being $S$ in the former, the argument of the former is a vector and the argument of the latter is a matrix). 
A function $f : \reals \to \reals^{n \times p}$ is ${\cal L}_\infty$ ($f \in {\cal L}_\infty$) if 
$\|f\|_\infty \coleq \esssup \{ |f(x)| : x \in \reals \}$ 
is finite.
%
%
%
The inverse of a quaternion $q \in \reals^4$ is denoted by $q^{-1}$ and the product of quaternions $q_1, q_2 \in \reals^4$ is denoted by $q_1 \otimes q_2$.

%



\subsection{Preliminaries on hybrid systems} \label{sec:COMHybridSystems}
In this paper, a hybrid system $\mathcal{H}$ is defined by $(C,F,D,G)$ as \cite{220}
\begin{equation} \label{eq:HyEq}
{\cal H} : 
\left\{
\begin{aligned}
    \dot \xi &= F(\xi) &\qquad \xi &\in C \\
    \xi^+ &= G(\xi)        &\qquad \xi &\in D
\end{aligned}
\right.
\end{equation}
where $\xi \in \reals^n$ is the state, $F : C \to \reals^n$ is the flow map defining a differential equation capturing the continuous dynamics, and $C \subset \reals^n$ defines the flow set on which flow is permitted. The mapping $G  : D \to \reals^n$ is the jump map defining the law resetting $\xi$ at jumps, and $D \subset \reals^n$ is the jump set on which jumps are permitted.

A solution $\xi$ to $\HS$ is a hybrid arc that is parameterized by $(t,j) \in \realsgeq \times \nats$, where $t$ is the elapsed ordinary time and $j$ is the number of jumps that have occurred. The domain of $\xi$, denoted by $\dom \xi \subset \realsgeq \times \nats$, is a {\it hybrid time domain}, in the sense that for every $(T,J) \in \dom \xi$, there exists a nondecreasing sequence $\{t_j\}_{j=0}^{J+1}$ with $t_0 = 0$ such that
\begin{equation*}
\dom \xi \cap \bigpar{[0,T] \times \{ 0, 1, \dots, J \}} = {\textstyle \bigcup\limits_{j = 0}^J} \bigpar{[t_j,t_{j+1}], j}.
\end{equation*}

A solution $\xi$ to $\HS$ is called {\it maximal} if it cannot be extended -- that is, if there does not exist another solution $\xi'$ to $\HS$ such that $\dom \xi$ is a proper subset of $\dom \xi'$ and $\xi(t,j) = \xi'(t,j)$ for all $(t,j) \in \dom \xi$. A solution is called {\it complete} if its domain is unbounded.

\section{Hybrid Finite-Time Estimation Algorithm} \label{sec:FTUniningControl}
\subsection{Problem Statement} \label{sec:FTProblemStatement}

In this paper, we develop a hybrid algorithm for estimating parameters of a continuous-time system of the form
\begin{align} \label{eqn:FTxdot}
    \dot x &= f(x,u(t)) + \phi(t) \theta 
\end{align}
where $x \in \mathbb{R}^n$ is the known state vector, $t \mapsto u(t) \in \reals^r$ is the known input, $t \mapsto \phi(t) \in \mathbb{R}^{n \times p}$ is the known regressor, $(x,u) \mapsto f(x,u) \in \reals^n$ is a known continuous function, $\theta \in \mathbb{R}^{p}$ is a vector of unknown constant parameters, and $n, r, p \in \nats$.

We propose to estimate $\theta$ using a hybrid algorithm, denoted by $\HS$ as in \eqref{eq:HyEq}, with data designed to solve the following problem.

\vspace{.05cm}
\noindent {\bf Problem Statement:} 
Given $t \mapsto \phi(t) \in \reals^{n \times p}$, $t \mapsto u(t) \in \reals^r$, and $(x,u) \mapsto f(x,u) \in \reals^n$, design the data $(C,F,D,G)$ of $\HS$ and determine conditions on $\phi$ that ensure the parameter estimate $\hat \theta$ converges in finite time to the unknown parameter vector $\theta$ in \eqref{eqn:FTxdot}.
\vspace{.05cm}

\subsection{Problem Solution} \label{sec:FTProblemSolution}
Given $t \mapsto \phi(t) \in \reals^{n \times p}$, $t \mapsto u(t) \in \reals^r$, and $(x,u) \mapsto f(x,u) \in \reals^n$, we design the data $(C, F, D, G)$ of $\HS$ so that, each time $\HS$ jumps, the parameter estimate $\hat \theta$ is equal to $\theta$ after the jump. 
To explain this construction, we begin by defining the dynamics of the state variables of $\HS$. 

During flows, we convert the plant dynamics \eqref{eqn:FTxdot} into a form similar to a linear regression model using state variables $\psi \in \reals^{n \times p}$ and $\eta \in \reals^{n}$, with dynamics \cite{Panteley2002Detectability}
\begin{equation} \label{eqn:FTtemp1}
\begin{aligned}
    \dot \psi &= -\lambda \psi + \phi(t) \\
    \dot \eta &= -\lambda (x + \eta) - f(x,u(t))
\end{aligned}
\end{equation}
where $\lambda > 0$ is a design parameter.
Defining
\begin{align} \label{eqn:FTeps}
    \varepsilon &\coleq x + \eta - \psi \theta, &
    y &\coleq x + \eta
\end{align}
it follows that $\psi$ and $\varepsilon$ are related via
\vspace{-2pt}
\begin{align} \label{eqn:FTregressionmodelC}
    y = \psi \theta + \varepsilon.
\end{align}
Since $\theta$ is constant, differentiating $\varepsilon$ along trajectories of \eqref{eqn:FTxdot}-\eqref{eqn:FTtemp1} yields
\vspace{-8pt}
\begin{align} \label{eqn:FTepsDot}
    \dot \varepsilon = - \lambda \varepsilon
\end{align}
%
%
%
Thus, $\varepsilon$ converges exponentially to zero during flows. 
Moreover, from \eqref{eqn:FTeps} and \eqref{eqn:FTregressionmodelC}, we have the following series of equivalences: 
%
$
    \varepsilon~\rightarrow~0
    \iff x + \eta~\rightarrow~\psi \theta
    \iff y~\rightarrow~\psi \theta.
$
%
Hence, $y$, $\psi$, and $\theta$ are related via a linear regression model plus an exponentially convergent term. We estimate $\theta$ using the classical gradient descent algorithm \cite{Narendra_adaptive_1989}, namely,
\begin{equation*} 
    \dot{\hat \theta} = \gamma \psi^\top (y - \psi \hat \theta) = \gamma \psi^\top (x + \eta - \psi \hat \theta)
\end{equation*}
where $\gamma > 0$ is a design parameter. Defining the parameter estimation error as $\tilde \theta \coleq \theta - \hat \theta$, we obtain
%
\vspace{-2pt}
\begin{equation} \label{eqn:FTtildethetaC}
    \dot {\tilde \theta} = - \gamma \psi^\top (\psi \tilde \theta + \varepsilon).
\end{equation}
Next, we study the finite-time convergence of the estimation error to zero. 
Let $(t,j) \mapsto \xi(t,j)$ be a solution to $\HS$ -- hence, defined on a hybrid time domain -- and consider the initial interval of flow $I^0 \coleq \{t : (t,0) \in \dom \xi \}$ with initial conditions 
$\psi(0,0) = 0$ and 
$\eta(0,0) = -x(0,0)$. 
From such initial conditions, we have from \eqref{eqn:FTtildethetaC} that, for all $t \in I^0$,
\vspace{-2pt}
\begin{align} \label{phi12statetransition}
    \tilde \theta(t,0) 
    &= \Phi(t,0) \tilde \theta(0,0) - \gamma \int_{0}^{t} \Phi(s,0) \psi^\top(s,0) \varepsilon(s,0) ds \notag \\
    &= \Phi(t,0) \tilde \theta(0,0)
\end{align}
where $(t,j) \mapsto \Phi(t,j)$ is the state transition matrix for $\tilde \theta$ generated by
%
$
    \dot{\Phi} = - \gamma \psi^\top \psi \Phi
$
 with 
$
    \Phi(0,0) = I, 
$
%
and the second equality in \eqref{phi12statetransition} follows from the fact that, since $\varepsilon(0,0) = x(0,0) + \eta(0,0) - \psi(0,0) \theta = 0$ and $\dot \varepsilon = - \lambda \varepsilon$, $\varepsilon(t,0) = 0$ for all $t \in I^0$. 
Then, if there exists a time $t_1 \in I^0$ such that the matrix $\Phi(t_1,0)-I$ is invertible, we reset $\hat \theta$ to the value of the function
\begin{equation} \label{eqn:FTR}
    R(\xi) \coleq K_1(\xi) \hat \theta + K_2(\xi) m
\end{equation}
where
\begin{align} \label{eqn:FTKi1}
    K_1(\xi) &\coleq - \big( \Phi - I \big)^{-1}, &
    K_2(\xi) &\coleq I - K_1(\xi)
\end{align}
and $m \in \reals^p$ is a memory state that stores the initial condition of $\hat \theta$. For readability, we define $\tilde m \coleq \theta - m$ and omit the argument $\xi$ of $K_1$ and $K_2$. Then, after the reset at hybrid time $(t_1,0)$, we obtain
\begin{align} \label{hatz_jump}
        &\hat \theta(t_1,1)
        = R(\xi(t_1,0))
        = K_1 \hat \theta(t_1,0) + K_2 m(t_1,0) \\
        &\quad= K_1 (\theta - \tilde \theta(t_1,0)) + K_2 (\theta - \tilde m(t_1,0)) \notag \\
        %
        %
        &\quad= - K_1 (\tilde \theta(t_1,0) - \tilde m(t_1,0)) - \tilde m(t_1,0) + (K_1 + K_2) \theta \notag \\
        &\quad= - K_1 (\Phi(t_1,0) - I) \tilde m(0,0) - \tilde m(0,0) + \theta \notag \\
        &\quad= \tilde m(0,0) - \tilde m(0,0) + \theta \notag
        = \theta \notag
\end{align} 
which follows from \eqref{phi12statetransition}-\eqref{eqn:FTKi1} and the fact that, since $m$ is constant during flows, $\tilde m(t_1,0) = \tilde m(0,0) = \tilde \theta(0,0)$. Hence, we have finite-time convergence of $\hat \theta$ to $\theta$.

In Theorem~\ref{thm:FTClosed_Loop_Stab_HS1} below, we establish that, if the integral of $\psi^\top \psi$ is positive definite, then there exists a time $t_1 \in I^0$ such that $\Phi(t_1,0)-I$ is invertible and hence $K_1$, $K_2$ in \eqref{eqn:FTKi1} are well defined. Thus, we introduce a state component, $Q$, with dynamics that compute the integral of $\psi^\top \psi$ online as
%
$
    \dot{Q} = \psi^\top \psi.
$
%

We define the jump map $G$ of $\HS$ to compute $\theta$ as in \eqref{hatz_jump} by resetting $\hat \theta$ to the value of $R$ in \eqref{eqn:FTR} at jumps. 
The state components $\psi$ and $\eta$ are reset as $\psi^+ = 0$ and $\eta^+ = -x$ which, from \eqref{eqn:FTeps}, implies that $\varepsilon^+ = x^+ + \eta^+ - \psi^+ \theta = 0$. 
The state component $\Phi$ is reset to the identity, which ensures that, during flows, $\Phi$ is the state transition matrix for the estimation error system \eqref{eqn:FTtildethetaC}.
The state component $Q$ is reset to zero, which ensures that $Q$ gives the value of the integral of $\psi^\top \psi$ over each interval of flow.
The state component $m$ is reset to the value of $R$ so that $m$ stores the value of $\hat \theta$ after each jump.
The flow and jump sets of $\HS$ are defined so that the algorithm jumps only when the minimum eigenvalue of $Q$ is greater than a design parameter, $\alpha > 0$, which ensures that the matrix $\Phi - I$ in \eqref{eqn:FTKi1} is invertible at jumps. 

We define the data of $\HS$ to model the closed-loop interconnection of the plant \eqref{eqn:FTxdot} and the proposed finite-time estimation algorithm. Thus, $\HS$ has state vector $\xi = (x, \hat \theta, \psi, \eta, \Phi, Q, m, \tau) \in {\cal X} \coleq \reals^n \times \reals^{p} \times \reals^{n \times p} \times \reals^{n} \times \reals^{p \times p} \times \reals^{p \times p} \times \reals^{p} \times \realsgeq$, where $x$ is the state of the plant, $\hat \theta$ is the estimate of $\theta$, and $\psi, \eta, \Phi, Q, m, \tau$ are auxiliary state variables.\footnote{The state component $\tau$ has dynamics such that it evolves as ordinary time $t$. Including $\tau$ in $\xi$ allows $\phi$ and $u$ to be part of the definitions of $F$ and $G$, rather than modeled as inputs to $\HS$. Thus, we can express $\HS$ as an autonomous hybrid system which allows us to leverage recent results on stability and robustness properties for such systems \cite{65,220}.} 
The data of $\HS$ is
\begin{equation} \label{eqn:FTH1}
\mbox{$\squeezespaces{0.6}
\begin{aligned}
    {\cal H} : 
    \left\{
    \begin{aligned}
    \dot \xi &= 
    \begin{bmatrix}
        f(x,u(\tau)) + \phi(\tau) \theta \\
        \gamma \psi^\top (x + \eta - \psi \hat \theta) \\ 
        -\lambda \psi + \phi(\tau) \\
        -\lambda (x + \eta) - f(x,u(\tau)) \\
        - \gamma \psi^\top \psi \Phi \\ 
        \psi^\top \psi \\
        0 \\
        1
    \end{bmatrix} =: F(\xi) &
    \xi &\in C \\
    \xi^+ &=
    (x, R(\xi), 0, -x, I, 0, R(\xi), \tau) 
    =: G(\xi) &
    \xi &\in D
    \end{aligned}
    \right.
\end{aligned}
$}
\end{equation}
with $R$ as in \eqref{eqn:FTR} and
\begin{equation} \label{eqn:FTCandD}
\begin{aligned} 
    C &\coleq \left\{ \xi \in {\cal X} : \lambda_{\min}(Q) \leq \alpha \right\} \\
    D &\coleq \left\{ \xi \in {\cal X} : \lambda_{\min}(Q) \geq \alpha \right\}
\end{aligned}
\end{equation}
%
%
%
where $\alpha > 0$ is a design parameter.

\begin{remark} \label{rem:FTConditionNumber}
The estimation algorithms proposed in \cite{hartman2012robust,242,256} compute a matrix inverse to achieve finite-time convergence of the parameter estimate to the true value. 
Similarly, the jump map $G$ of $\HS$ in \eqref{eqn:FTH1} requires computation of the inverse of $\Phi - I$ to evaluate the function $R$ -- see \eqref{eqn:FTR} and \eqref{eqn:FTKi1}.
However, in contrast to $\HS$, the algorithms in \cite{hartman2012robust,242,256} use the determinant to evaluate whether the matrix inverse is well defined, which may result in the matrix being numerically ill-conditioned for inversion.
For example, given ${d} > 0$, consider the matrix
\begin{align*}
    M &\coleq \frac{\sqrt{d\phantom{.}}}{3a - a^2} \matt{1+(1-a)^2 & 4-3a \\ 4-3a & 4+(2-a)^2}
\end{align*}
where $a \in (0, 3)$. 
For any such $a$, the determinant of $M$ is equal to ${d}$, but the condition number of $M$ grows unbounded as $a$ approaches $0$ or $3$. Recall that the condition number measures how sensitive matrix inversion is to numerical errors. This fact motivates using a method other than the determinant to evaluate whether a matrix is invertible. For $\HS$ in \eqref{eqn:FTH1}-\eqref{eqn:FTCandD}, we impose a lower bound on the minimum eigenvalue of the solution component $Q$ in the jump set, which allows us to upper-bound the condition number of $\Phi-I$ at jumps -- see Theorem~\ref{thm:FTClosed_Loop_Stab_HS1}.
%
%
\end{remark}

\begin{example} \label{exp:FTNumericalExample1}
To demonstrate how the proposed algorithm works, consider the following system inspired by \cite{Adetola_finite_time_est_2008}, with state $x = (x_1, x_2, x_3) \in \reals^3$ and dynamics 
\begin{equation} \label{eqn:FTexamplePlant}
\begin{aligned}
    \dot x_1 &= x_2 + \theta_1 x_1 \\
    \dot x_2 &= x_3 + \theta_2 x_1 \\
    \dot x_3 &= \theta_3 x_1^3 + \theta_4 x_2 + \theta_5 x_3 + (2 + \sin(x_1)) \, \kappa(x, u)
\end{aligned}
\end{equation}
where $\theta = (\theta_1, \theta_2, \cdots, \theta_5) \in \reals^5$ is a vector of unknown constant parameters,
the input is $u \coleq (x_r, \hat \theta)$ where $x_r \in \reals$ is a constant reference value and $\hat \theta \in \reals^5$ is an estimate of $\theta$,
and the control signal $(x, u) \mapsto \kappa(x, u) \in \reals$ is designed as in \cite{Jung-Shan1999Nonlinearities} so that the state component $x_1$ converges asymptotically to $x_r$.

To estimate $\theta$ using $\HS$, we express the dynamics of \eqref{eqn:FTexamplePlant} in the form of \eqref{eqn:FTxdot} with
\begin{equation} \label{eqn:FTexamplePhi}
\begin{aligned}
    f(x,u) &\coleq \matt{ x_2 \\ x_3 \\ \big( 2 + \sin(x_1) \big) \, \kappa(x, u) }
\end{aligned}
\end{equation}
and, for each solution pair $(x, u)$ to \eqref{eqn:FTexamplePlant},
\begin{equation} \label{eqn:FTexamplePhi2}
\begin{aligned}
    \phi(t) &\coleq \matt{x_1(t) & 0 & 0 & 0 & 0 \\
                        0 & x_1(t) & 0 & 0 & 0 \\
                        0 & 0 & x_1(t)^3 & x_2(t) & x_3(t)}
\end{aligned}
\end{equation}
for all $t \in \dom (x, u)$. 
We simulate $\HS$ in \eqref{eqn:FTH1} with $f$ and $\phi$ in \eqref{eqn:FTexamplePhi}-\eqref{eqn:FTexamplePhi2}. 
The system \eqref{eqn:FTexamplePlant} has parameters $\theta = (-1, -2, 1, 2, 3)$, reference value $x_r = 1$, and the initial condition for $x$ is $(0,0,0)$. 
The estimator has design parameters $\gamma = 1$, $\lambda = 0.1$, and $\alpha = 0.2$, and initial conditions $\hat \theta = m = (0,0,0,0,0)$, $\psi = 0$, $\eta = -x$, $\Phi = I$, and $Q = 0$. 
Verifying PE as in \eqref{eqn:FTcontinuousPE} for the regressor $\phi$ in \eqref{eqn:FTexamplePhi2} is difficult in practice since the property must hold for all future time. However, with these parameters and initial conditions, it can be shown numerically\footnote{Code at \href{https://github.com/HybridSystemsLab/HybridFT_NonlinearSystem}{\url{https://github.com/HybridSystemsLab/HybridFT_NonlinearSystem}}} that $\phi$ is such that $\lambda_{\min}(Q(2.93, 0)) \geq \alpha $. Hence, $\HS$ jumps and the parameter estimate $\hat \theta$ converges in finite-time to $\theta$ according to \eqref{hatz_jump}, as shown in the top subplot of Figure~\ref{fig:FTexampleEstimationError1}.
Following convergence of the parameter estimate, the state component $x_1$ converges to the reference value $x_r$ as shown in the bottom sub-plot of Figure~\ref{fig:FTexampleEstimationError1}. 
\begin{figure}[!htb]
    \centering
        \includegraphics[width=1\linewidth,height=\textheight,keepaspectratio]{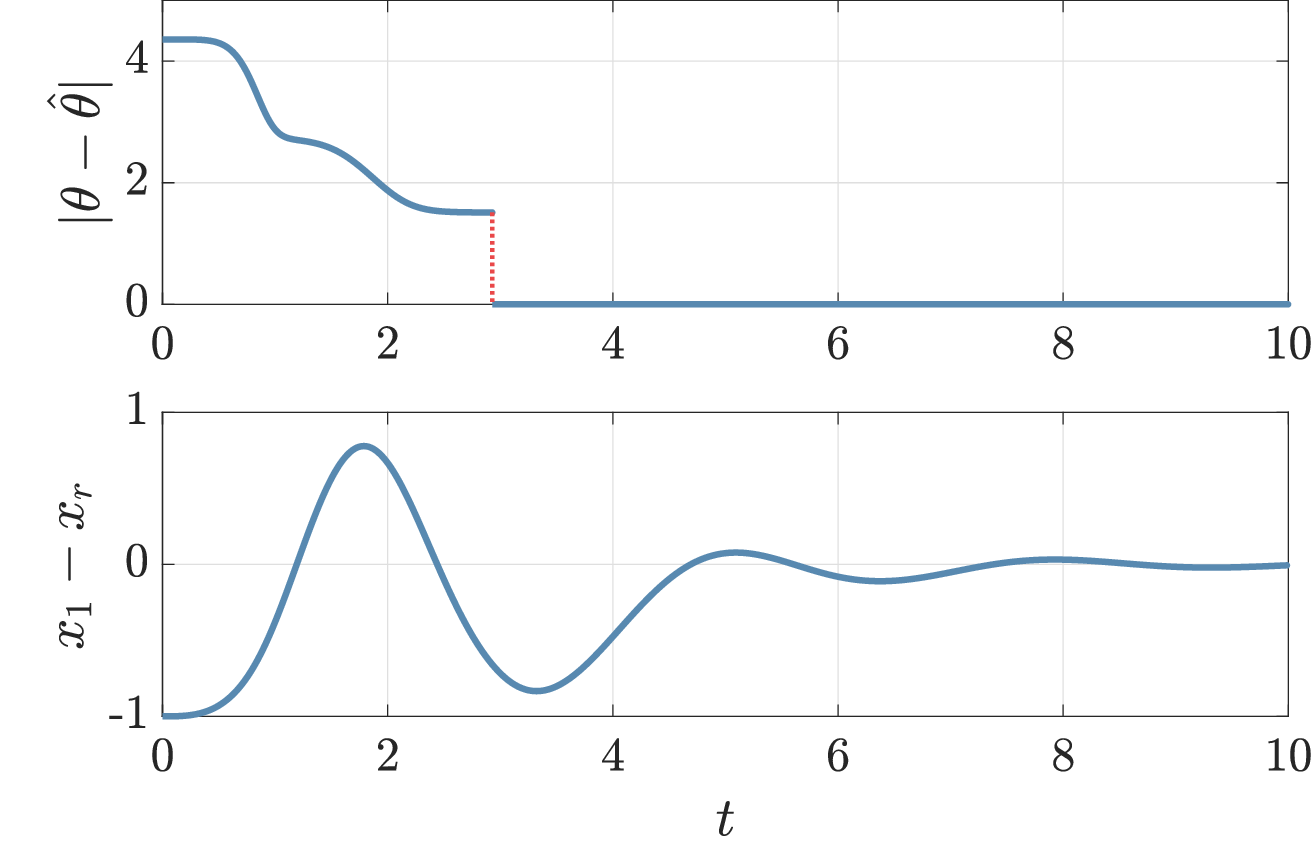}
        \caption{Top: The projection onto $t$ of the norm of the parameter estimation error for $\HS$. 
        Bottom: The error between the state component $x_1$ and the reference value $x_r$.}
    \label{fig:FTexampleEstimationError1}
\end{figure}
\end{example}
\vspace{-8pt} 

\section{Convergence Analysis} \label{sec:FTConvergenceAnalysis}

%

We now establish our main convergence result stating conditions ensuring that each maximal solution $\xi$ to $\HS$ with initial conditions in a closed set converges in finite time to the set
\begin{equation} \label{eqn:FTCalA}
    \A \coleq \defset{\xi \in \X}{\hat \theta = \theta}.
\end{equation}
We impose the following boundedness condition on the matrix signal $\phi$. This condition is analogous to the bound on the Euclidean norm that is typically imposed when analyzing the excitation properties of vector signals \cite{Narendra_adaptive_1989,Sastry1989Adaptive,Khalil2002Nonlinear}.
\begin{assumption} \label{asm:FTphiBound}
Given $t \mapsto \phi(t) \in \reals^{n \times p}$,
there exists $\phi_M > 0$ such that $|\phi(t)|_{\mathrm F} \leq \phi_M$ for all $t \geq 0$.\footnote{As written, Assumption~\ref{asm:FTphiBound} imposes completeness of $t~\mapsto~\phi(t)$. The assumption of completeness of $\phi$ can be omitted at the expense of increased complexity in the statements and proofs of the results herein.}
\end{assumption}
We impose the following finite excitation condition on $\phi$. Remarkably, as we will establish in Theorem~\ref{thm:FTClosed_Loop_Stab_HS1} below, this condition is sufficient to guarantee finite-time convergence of the parameter estimate to the true value for the proposed estimation algorithm.
\begin{assumption} \label{def:FTFE}
Given $t \mapsto \phi(t) \in \reals^{n \times p}$, there exist $t' \geq 0$, $T, \mu_\phi > 0$, and $L_\phi \geq 0$ such that
\begin{equation}  \label{eqn:FTFE}
    \int_{t'}^{t'+T} 
    \phi (s)^\top 
    \phi (s) ds
    \geq \mu_\phi I
\end{equation}
%
%
and, for all $s, v \in [t', t' + T]$,
\begin{equation}  \label{eqn:Lischitz_phi}
    |\phi(s) - \phi(v)| \leq L_\phi |s - v|.
\end{equation}
\end{assumption}
%
%
%
\begin{remark} \label{rem:FTExampleAssmp}
To demonstrate the practicality of Assumptions~\ref{asm:FTphiBound} and~\ref{def:FTFE}, recall Example~\ref{exp:FTNumericalExample1}. Verifying PE as in \eqref{eqn:FTcontinuousPE} for the regressor $\phi$ in \eqref{eqn:FTexamplePhi2} is difficult in practice since the property must hold for all future time. However, it can be shown numerically that $\phi$ satisfies Assumption~\ref{asm:FTphiBound} with $\phi_M = 6.642$, and satisfies Assumption~\ref{def:FTFE} with $t' = 0$ sec, $T = 2.93$ sec, $\mu_\phi = 1.85$, and $L_\phi = 15.9$.
\end{remark}
%
%
%
Before stating our main convergence result, we establish the following lemma that relates excitation properties of $\phi$ to those of $\psi$.
\begin{lemma} \label{lem:FT_FEPsi}
Given $t \mapsto \phi(t) \in \reals^{n \times p}$, suppose that Assumptions~\ref{asm:FTphiBound} and~\ref{def:FTFE} hold with $\phi_M,\mu_\phi > 0$, $t'\geq 0$, $T > 0$, and $L_\phi \geq 0$. Then, for each $\lambda > 0$, the solution $\psi : [t', \infty) \to \reals^{n\times p}$ to the system
\begin{equation} \label{eqn:FTpsidynamics0}
\begin{aligned} 
    \dot \psi = -\lambda \psi + \phi(t), \qquad \psi(t') = 0
\end{aligned}
\end{equation}
satisfies
\begin{equation}  \label{eqn:FT_psiFE}
    \int_{t'}^{t'+T} 
    \psi (s)^\top 
    \psi (s) ds
    \geq \mu_\psi I
\end{equation}
%
%
where
\begin{equation} \label{eqn:FT_mu_a}
\!\!\!
\mbox{$\squeezespaces{0.4}
\begin{aligned}
\mu_\psi \coleq \frac{\mu_\phi}{4(2a - \lambda)^2}, \quad \!\!\!
a \coleq \max \left\{ \lambda, \frac{16 \phi_M^2}{\mu_\phi}, 4 L_\phi \sqrt{ \frac{T}{\mu_\phi} } \right\}.
\end{aligned}
$}
\end{equation}
\end{lemma}
\begin{proof}
    This proof is in Appendix~\ref{apx:FT_FEPsi}.
\end{proof}
%
Given the solution $\psi$ to \eqref{eqn:FTpsidynamics0}, we define
\begin{equation} \label{eqn:FT_mubar}
    \mu \coleq \lambda_{\min} \left( \int_{t'}^{t'+T} \psi(s)^\top \psi(s) ds \right).
\end{equation}
The value of $\mu$ is determined only by $t \mapsto \phi(t)$, $\lambda$, $t'$, and $T$ and, in particular, does not depend on solutions to $\HS$.
Assumptions~\ref{asm:FTphiBound} and~\ref{def:FTFE} allow us to establish a positive lower bound on $\mu$ that can be computed a priori: by Lemma~\ref{lem:FT_FEPsi}, $\mu~\geq~\mu_\psi~>~0$ with $\mu_\psi$ as in \eqref{eqn:FT_mu_a}.
Hence, the design parameter $\alpha$ in \eqref{eqn:FTCandD} can always be chosen such that $0~<~\alpha~<~\mu$, which, as we will establish in Theorem~\ref{thm:FTClosed_Loop_Stab_HS1} below, guarantees finite-time convergence of the parameter estimate to the true value for the proposed estimation algorithm.
\begin{theorem} \label{thm:FTClosed_Loop_Stab_HS1}
Given $t \mapsto \phi(t) \in \reals^{n \times p}$, $t \mapsto u(t) \in \reals^r$, $(x,u) \mapsto f(x,u) \in \reals^n$ in \eqref{eqn:FTxdot}, 
and $\lambda, \alpha > 0$,
suppose that Assumptions~\ref{asm:FTphiBound} and~\ref{def:FTFE} hold with $t' = 0$, $T, \phi_M, \mu_\phi > 0$, $L_\phi \geq 0$, and $\alpha < \mu$ with $\mu$ as in \eqref{eqn:FT_mubar}.
%
%
Then, for each $\gamma > 0$, the hybrid system $\HS$ in \eqref{eqn:FTH1} is such that the parameter estimation error $(t,j) \mapsto \tilde \theta(t,j) \coleq \theta - \hat \theta(t,j)$ for each solution $\xi$ to $\HS$ from $\xi(0,0) \in {\cal X}_0 \coleq \{\xi \in {\cal X} : \hat \theta = m, \, \psi = 0, \, \eta = -x, \, \Phi = I, \, Q = 0, \, \tau = 0 \}$ satisfies
\begin{equation} \label{eqn:FTxibound}
    |\tilde \theta(t,j)| \leq |\tilde \theta(0,0)| \quad \forall (t,j) \in \dom \xi
\end{equation}
%
%
and the parameter estimate satisfies
\begin{equation} \label{eqn:FTFTbound}
   \hat \theta(t,j) = \theta
\end{equation}
for all $(t,j) \in \dom \xi$ satisfying $t \geq T$ and $j \geq 1$. 
Moreover, for each maximal solution $\xi$ to $\HS$, the condition number of $\Phi - I$ immediately before the first jump at hybrid time $(t_1,0)$ satisfies
\begin{equation} \label{eqn:FTdeltabound}
    \text{\rm cond}(\Phi(t_1,0) - I) \leq \frac{1 + \sigma(\gamma,\alpha)}{1 - \sigma(\gamma,\alpha)}
\end{equation}
where $t_1 \in (0, T]$ and
\begin{align} \label{eqn:FTalphabeta}
    \sigma(\gamma,\alpha) &\coleq \sqrt{ 1 - \frac{2 \gamma \alpha}{(1 + \gamma \psi_M^2 T)^2}}, & 
    \psi_M &\coleq \frac{\phi_M}{\lambda}.
\end{align}
%
%
%
%
\end{theorem}
\begin{remark}
    In practice, the integral in \eqref{eqn:FT_mubar} does not need to be computed a priori to implement the proposed estimation scheme.
    This is due to the fact that, for each solution $\xi$ to $\HS$ from $\xi(0,0) \in {\cal X}_0$, the solution component $\psi$ is the solution to \eqref{eqn:FTpsidynamics0} with $t' = 0$ on $I^0 \coleq \{t : (t,0) \in \dom \xi\}$.
    Thus, the solution component $Q$ satisfies $Q(t,0) = \int_{0}^{t} \psi(s,0)^\top \psi(s,0) ds$ for all $t \in I^0$, and the jump set $D$ in \eqref{eqn:FTCandD} triggers a jump only when $Q(t,0) \geq \alpha I$.
    Hence, $\HS$ jumps opportunistically when reconstruction of $\theta$ is possible through the conditions imposed in Assumptions~\ref{asm:FTphiBound} and~\ref{def:FTFE}.
\end{remark}
\begin{example}
Recall Example~\ref{exp:FTNumericalExample1}.
It can be shown numerically that $\mu$ in \eqref{eqn:FT_mubar} is greater than $0.2$, and thus $\alpha < \mu$, whereas the bound in \eqref{eqn:FT_mu_a} is 
$\mu_\psi = 7.9 \times 10^{-7}$.
Hence, from Remark~\ref{rem:FTExampleAssmp} it follows that the conditions of Theorem~\ref{thm:FTClosed_Loop_Stab_HS1} are satisfied and, in accordance with Theorem~\ref{thm:FTClosed_Loop_Stab_HS1}, the parameter estimate $\hat \theta$ converges in finite time to $\theta$ as shown in the top subplot of Figure~\ref{fig:FTexampleEstimationError2}. The bottom subplot of Figure~\ref{fig:FTexampleEstimationError2} shows the condition number of $\Phi - I$ in blue along with the bound \eqref{eqn:FTdeltabound} in dashed black. 
In this case, \eqref{eqn:FTdeltabound} ensures that cond$(\Phi(t_1,0) - I) \leq 1.67 \times 10^9$ and the actual value is cond$(\Phi(t_1,0) - I) = 26.4$.
\begin{figure}[!htb]
    \centering
        \includegraphics[width=1\linewidth,height=\textheight,keepaspectratio]{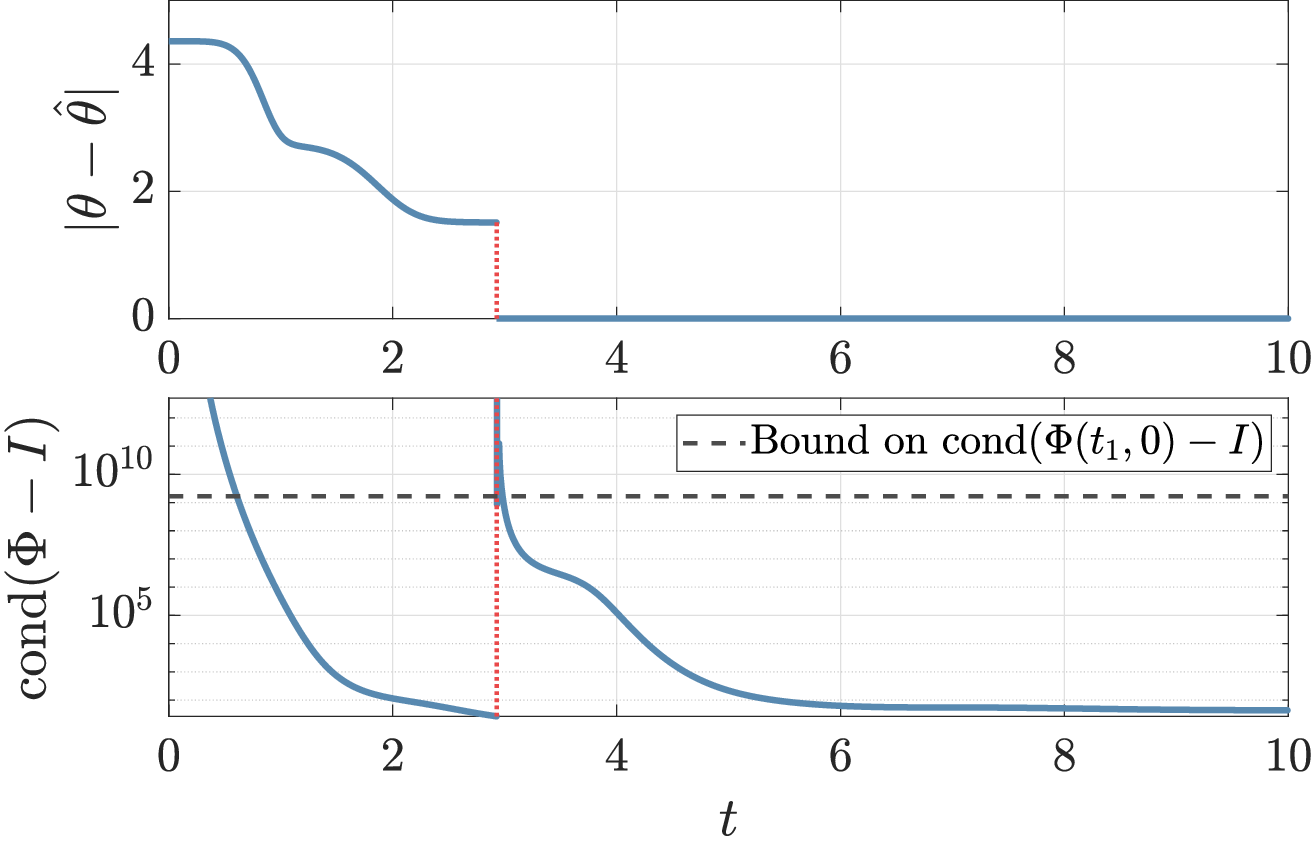}
        \caption{Top: The projection onto $t$ of the norm of the parameter estimation error for $\HS$. 
        Bottom: The projection onto $t$ of cond$(\Phi - I)$ in blue along with the upper bound on cond$(\Phi(t_1,0) - I)$ from \eqref{eqn:FTdeltabound} in dashed black.}
    \label{fig:FTexampleEstimationError2}
\end{figure}
\end{example}
\subsection{Proof of Theorem~\ref{thm:FTClosed_Loop_Stab_HS1}}

The proof of Theorem~\ref{thm:FTClosed_Loop_Stab_HS1} proceeds as follows. 
First, in Lemma~\ref{lem:FTPhiExcitation}, we establish a uniform upper bound on the norm of a continuous-time signal $t \mapsto \psi(t)$ that shares the same dynamics as the flow map for the state component $\psi$ in \eqref{eqn:FTH1}. 
Next, in Lemma~\ref{lem:FTPhiStability}, we establish a uniform upper bound on the norm and on the condition number of a continuous-time signal $t \mapsto \Omega(t)$ that shares the same dynamics as the flow map for the state component $\Phi$ in \eqref{eqn:FTH1}. 
Finally, these results are used to prove Theorem~\ref{thm:FTClosed_Loop_Stab_HS1}.

%
\begin{lemma} \label{lem:FTPhiExcitation}
Given $t \mapsto \phi(t) \in \reals^{n \times p}$, suppose that Assumption~\ref{asm:FTphiBound} holds with $\phi_M > 0$. 
Then, for each $t' \geq 0$ and each $\lambda > 0$, the solution $\psi : [t',\infty) \to \reals^{n \times p}$ to \eqref{eqn:FTpsidynamics0} satisfies
\begin{equation} \label{eqn:FTpsibound}
    |\psi(t)| \leq \frac{\phi_M}{\lambda} =: \psi_M
    \qquad \forall t \geq t'.
\end{equation}
\end{lemma}
\begin{proof}
%
%
Consider the Lyapunov function
\begin{equation} \label{eqn:FTVpsi}
    W(\psi) 
    \coleq \frac{1}{2} \tr (\psi^\top \psi) = \frac{1}{2} |\psi|_{\mathrm F}^2
    \quad \forall \psi \in \reals^{n \times p}.
\end{equation}
Differentiating $W$ along solutions to \eqref{eqn:FTpsidynamics0} , we obtain
\begin{align*} %
    \dot W(\psi)
    &= \frac{1}{2} \tr (\dot \psi^\top \psi + \psi^\top \dot \psi)
    = \tr (- \lambda \psi^\top \psi + \psi^\top \phi(t)) \\
    &= - 2 \lambda W(\psi) + \tr(\psi^\top \phi(t))
\end{align*}
where the second equality follows from the flow map for $\psi$ in \eqref{eqn:FTpsidynamics0} , and the third equality follows from the definition of $W$ in \eqref{eqn:FTVpsi}. 
Applying the Cauchy-Schwarz inequality yields
\begin{equation*} %
    \dot W(\psi)
    \leq - 2 \lambda W(\psi) + \sqrt{2 W(\psi)} |\phi(t)|_{\mathrm F}.
\end{equation*}
Hence,
\begin{equation*} 
    \dot W(\psi) 
    \leq 0 \ \
    \forall \psi \in \reals^{n \times p} : 
    W(\psi) \geq \frac{1}{2 \lambda^2} |\phi(t)|_{\mathrm F}^2.
\end{equation*}
Let $\psi$ be the solution to \eqref{eqn:FTpsidynamics0}. Using the bounds above, we conclude that
\begin{equation*}
    W(\psi(t))  
    \leq W(\psi(t')) + \frac{1}{2 \lambda^2} \phi_M^2
    = \frac{1}{2 \lambda^2} \phi_M^2
\end{equation*}
for all $t \geq t'$, where the first inequality follows from Assumption~\ref{asm:FTphiBound}, and the last equality follows from the fact that $\psi(t') = 0$. Using the definition of $W$ in \eqref{eqn:FTVpsi}, we obtain
\begin{equation*}
    |\psi(t)| \leq 
    |\psi(t)|_F 
    = \sqrt{2 W(\psi(t))}
    \leq \sqrt{2 \frac{\phi_M^2}{2 \lambda^2}}
    = \frac{\phi_M}{\lambda}
\end{equation*}
for all $t \geq t'$. Hence, \eqref{eqn:FTpsibound} holds.
\end{proof}
\begin{lemma} \label{lem:FTPhiStability}
Given $t \mapsto \phi(t) \in \reals^{n \times p}$, $\lambda > 0$, and $t' \geq 0$, let $\psi : [t' , \infty) \to \reals^{n \times p}$ be the solution to \eqref{eqn:FTpsidynamics0}, and consider the system
\begin{align} \label{eqn:FTlemsysa}
    \dot{\vartheta} &= - \gamma \psi(t)^\top \psi(t) \vartheta
\end{align}
where $\gamma > 0$. 
Let $\Omega : \realsgeq \times \realsgeq \to \reals^{p \times p}$ be the state transition matrix of \eqref{eqn:FTlemsysa}. Then, for each $t \geq 0$,
\begin{equation} \label{eqn:FTOmegaUpperBound}
    |\Omega(t'+t,t')| \leq 1.
\end{equation}
Moreover, given $\alpha > 0$, suppose that Assumptions~\ref{asm:FTphiBound} and~\ref{def:FTFE} hold with $T, \phi_M, \mu_\phi > 0$, $L_\phi \geq 0$, and $\alpha < \mu$ with $\mu$ as in \eqref{eqn:FT_mubar}. Then,
\begin{equation} \label{eqn:FTOmegaTBound}
    |\Omega(t'+T,t')| \leq \sigma(\gamma,\alpha) < 1
\end{equation}
where
\begin{equation} \label{eqn:FTsigmaLemma}
    \sigma(\gamma,\alpha) \coleq \sqrt{1 - \frac{2 \gamma \alpha}{(1 + \gamma \psi_M^2 T)^2}}
\end{equation}
with $\psi_M$ as in \eqref{eqn:FTpsibound}. Additionally, the matrix
\iftoggle{report}{}{\linebreak}
$\Omega(t'+T,t') - I$ is invertible with
\begin{equation} \label{eqn:FTK1Bound}
    |(\Omega(t'+T,t') - I)^{-1}| 
    \leq \frac{1}{1 - \sigma(\gamma,\alpha)}
\end{equation}
and the condition number of $\Omega(t'+T,t') - I$ satisfies
\begin{equation} \label{eqn:FTcondOmega}
    \text{\rm cond}(\Omega(t'+T,t') - I) \leq \frac{1 + \sigma(\gamma,\alpha)}{1 - \sigma(\gamma,\alpha)}.
\end{equation}
\end{lemma}
\begin{remark}
The property $\sigma(\gamma,\alpha) < 1$ in \eqref{eqn:FTOmegaTBound} follows from the fact that, since $\alpha < \mu$ with $\mu$ as in \eqref{eqn:FT_mubar}, $\alpha I \leq \mu I \leq \int_{t'}^{t'+T} \psi(s)^\top \psi(s) ds \leq \psi_M^2 T I$, with $\psi_M$ as in \eqref{eqn:FTpsibound}. Thus,
\begin{equation*}
    0 < \frac{2 \gamma \alpha}{(1 + \gamma \psi_M^2 T)^2}
    \leq \frac{2 \gamma \psi_M^2 T}{(1 + \gamma \psi_M^2 T)^2}
    \leq \frac{1}{2}.
\end{equation*}
\end{remark}
To prove Lemma~\ref{lem:FTPhiStability}, we first recall the following result from \cite{Horn2012Matrix}.
\begin{lemma} \label{lem:FTInvIminusA}
Given a matrix $A \in \reals^{p \times p}$, if $|A| < 1$ then $I - A$ is invertible and $|(I - A)^{-1}| \leq (1 - |A|)^{-1}$.
\end{lemma}
\begin{proof}[\unskip\nopunct]{\it Proof of Lemma~\ref{lem:FTPhiStability}:}
%
%
Consider the following Lyapunov function:
\begin{equation} \label{eqn:FTV1}
    V(\vartheta) \coleq \frac{1}{2 \gamma} \vartheta^\top \vartheta
    \quad \forall \vartheta \in \reals^{p}.
\end{equation}
We first show that \eqref{eqn:FTOmegaUpperBound} holds. Differentiating $V$ along trajectories of \eqref{eqn:FTlemsysa}, we obtain
\begin{equation} \label{eqn:FTVdot}
    \dot V(\vartheta) 
    = -\vartheta^\top \psi^\top \psi \vartheta
    = -\big| \psi \vartheta \big|^2
    \leq 0.
\end{equation}
Let $\vartheta$ be the solution to \eqref{eqn:FTlemsysa}. Using the bound above,
\begin{equation*}
    V(\vartheta(t'+t)) \leq V(\vartheta(t'))
    \quad \forall t \geq 0.
\end{equation*}
Using the definition of $V$ in \eqref{eqn:FTV1}, we obtain that $|\vartheta(t'+t)| \leq |\vartheta(t')|$ for all $t \geq 0$, which we rewrite as
\begin{equation*}
    |\Omega(t'+t,t') \vartheta(t')| \leq |\vartheta(t')|
    \quad \forall t \geq 0
\end{equation*}
which, if $|\vartheta(t')| \neq 0$, implies that
\begin{equation*}
    \frac{|\Omega(t'+t,t') \vartheta(t')|}{|\vartheta(t')|} \leq 1.
\end{equation*}
Since this inequality holds for any $\vartheta(t') \in \reals^{p} \setminus \{0\}$, it follows that
\begin{equation*}
    |\Omega(t'+t,t')| = \sup_{s \in \reals^{p} \setminus \{0\}} 
    \frac{|\Omega(t'+t,t') s|}{|s|} \leq 1
\end{equation*}
for all $t \geq 0$. Hence, \eqref{eqn:FTOmegaUpperBound} holds.

To show that \eqref{eqn:FTOmegaTBound} holds, we follow an approach inspired by \cite{242}. Integrating \eqref{eqn:FTVdot} along the solution $\vartheta$ to \eqref{eqn:FTlemsysa} between times $t'$ and $t' + T$, we obtain
\begin{align*}
&V(\vartheta(t'+T)) - V(\vartheta(t'))
= - \int_{t'}^{t'+T} \left| \psi(s) \vartheta(s) \right|^2 ds \\
& \ = \mbox{$\squeezespaces{0.5}- \displaystyle\int_{t'}^{t'+T} \left| \psi(s) \left( \vartheta(t') - \gamma \int_{t'}^{s} \psi^\top(u) \psi(u) \vartheta(u) du \right) \right|^2 ds $} 
%
%
%
%
\end{align*}
where the second equality follows from the fact that the solution $t \mapsto \vartheta(t)$ to \eqref{eqn:FTlemsysa} is 
$\vartheta(t) = \vartheta(t') - \gamma \int_{t'}^{t} \psi^\top(u) \psi(u) \vartheta(u) du$ for all $t \geq t' \geq 0$.
Using the fact that, for any $a, b \in \reals$ and any $\varrho > 0$, $(a-b)^2 \geq \frac{\varrho}{1+\varrho}a^2 - \varrho b^2$, we have
\begin{align*}
&V(\vartheta(t'+T)) - V(\vartheta(t')) \\
%
%
&\quad\leq - \frac{\varrho}{1 + \varrho} \int_{t'}^{t'+T} \left| \psi(s) \vartheta(t') \right|^2 ds \\
&\qquad + \varrho \int_{t'}^{t'+T} \left| \gamma \psi(s) \int_{t'}^{s} \psi^\top(u) \psi(u) \vartheta(u) du \right|^2 ds. 
%
%
%
\end{align*}
Using Lemma~\ref{lem:FTPhiExcitation} and the triangle and Cauchy-Schwarz inequalities, we have
\begin{align*}
V(\vartheta(t'+T&)) - V(\vartheta(t')) \\
%
%
%
\leq &- \frac{\varrho}{1 + \varrho} \int_{t'}^{t'+T} \left| \psi(s) \vartheta(t') \right|^2 ds \\
&+ \varrho \gamma^2 \psi_M^4 \int_{t'}^{t'+T} T \int_{t'}^{t'+T} \left| \psi(u) \vartheta(u) \right|^2 du ds. 
%
%
\end{align*}
To bound the first term on the right-hand-side of the expression above, we use the fact that, from \eqref{eqn:FT_mubar} and since $\alpha < \mu$, 
\begin{align*}
&\int_{t'}^{t'+T} |\psi(s) \vartheta(t')|^2 ds \\
&\qquad = \vartheta(t') \bigg( \int_{t'}^{t'+T} \psi(s)^\top \psi(s) ds \bigg) \vartheta(t') \geq \alpha |\vartheta(t')|^2.
\end{align*}
Thus,
\begin{align*}
&V(\vartheta(t'+T)) - V(\vartheta(t')) \\
%
%
%
%
&\ \leq - \frac{\varrho \alpha}{1 + \varrho} \big| \vartheta(t') \big|^2 + \varrho \gamma^2 \psi_M^4 T^2 \int_{t'}^{t'+T} \left| \psi(u) \vartheta(u) \right|^2 du \\
&\ \leq \mbox{$\squeezespaces{0.3} - \displaystyle\frac{2 \gamma \varrho \alpha}{1 + \varrho} V(\vartheta(t')) - \varrho \gamma^2 \psi_M^4 T^2 \big( V(\vartheta(t'+T)) - V(\vartheta(t')) \big) $}
\end{align*}
where the second inequality follows from the definition of $V$ in \eqref{eqn:FTV1}. 
Hence, we have
\begin{equation*} 
    V(\vartheta(t'+T)) \leq (1 - \beta) V(\vartheta(t'))
\end{equation*}
where 
\vspace{-2pt}
\begin{equation*}
    \beta \coleq \frac{2 \varrho \gamma \alpha}{(1 + \varrho)(1 + \varrho \gamma^2 \psi_M^4 T^2)}.
\end{equation*}
By choosing $\varrho = \frac{1}{\gamma \psi_M^2 T}$, we obtain
\begin{equation*}
    \beta = \frac{2 \gamma \alpha}{(1 + \gamma \psi_M^2 T)^2}
\end{equation*}
and, from \eqref{eqn:FT_mubar} and Lemma~\ref{lem:FTPhiExcitation} and since $\alpha < \mu$, it follows that 
\begin{equation*}
    \alpha I < \mu I \leq \int_{t'}^{t'+T} \psi(s)^\top \psi(s) ds \leq \psi_M^2 T I.
\end{equation*}
Hence, $\beta \in (0, 1)$. 
Using the definition of $V$ in \eqref{eqn:FTV1}, it follows that $|\vartheta(t'+T)| \leq \sqrt{1 - \beta} |\vartheta(t')|$, which we rewrite as
\begin{equation*}
    |\Omega(t'+T,t') \vartheta(t')| \leq \sqrt{1 - \beta} |\vartheta(t')|
\end{equation*}
which, if $|\vartheta(t')| \neq 0$, implies that
\begin{equation*}
    \frac{|\Omega(t'+T,t') \vartheta(t')|}{|\vartheta(t')|} \leq \sqrt{1 - \beta}.
\end{equation*}
Since this inequality holds for any $\vartheta(t') \in \reals^{p} \setminus \{0\}$, it follows that
\begin{equation*}
    |\Omega(t'+T,t')| = \sup_{s \in \reals^{p} \setminus \{0\}} 
    \frac{|\Omega(t'+T,t') s|}{|s|} \leq \sqrt{1 - \beta} < 1
\end{equation*}
where the last inequality follows from the fact that $\beta \in (0, 1)$. 
Hence, \eqref{eqn:FTOmegaTBound} holds and, in view of Lemma \ref{lem:FTInvIminusA}, the matrix $I - \Omega(t'+T,t')$ is invertible. Thus, $\Omega(t'+T,t') - I$ is invertible.

To upper bound the condition number of $\Omega(t'+T,t') - I$, we first use the triangle inequality to obtain
\begin{equation} \label{eqn:FTOmegaI}
    |\Omega(t'+T,t') - I| 
    \leq |\Omega(t'+T,t')| + |I|
    \leq \sigma(\gamma,\alpha) + 1.
\end{equation}
Next, since $|\Omega(t'+T,t')| \leq \sigma(\gamma,\alpha) < 1$, it follows from Lemma~\ref{lem:FTInvIminusA} that 
\begin{align} \label{eqn:FTOmegaIinv}
    &|(\Omega(t'+T,t') - I)^{-1}| \\
    &\qquad \qquad \qquad \leq \frac{1}{1 - |\Omega(t'+T,t')|}
    \leq \frac{1}{1 - \sigma(\gamma,\alpha)}. \notag
\end{align}
Hence, \eqref{eqn:FTK1Bound} holds.
Using the expressions in \eqref{eqn:FTOmegaI} and \eqref{eqn:FTOmegaIinv}, we conclude that \eqref{eqn:FTcondOmega} holds.
\end{proof}
\begin{lemma} \label{prop:FTepsStab}
Given the hybrid system $\HS$ in \eqref{eqn:FTH1}, for each $\lambda, \gamma > 0$ and each solution $\xi$ to $\HS$ from $\xi(0,0) \in {\cal X}_0$, $(t,j)~\mapsto~\varepsilon(t,j)$ in \eqref{eqn:FTeps} satisfies
\begin{equation} \label{eqn:FTepsbound}
    \varepsilon(t,j) = 0
    \quad \forall (t,j) \in \dom \xi.
\end{equation}
\end{lemma}
\begin{proof}
The proof of  Lemma~\ref{prop:FTepsStab} is immediate since, for each solution $\xi$ to $\HS$ from $\xi(0,0) \in {\cal X}_0$, $\varepsilon(0,0) = 0$ and, from \eqref{eqn:FTH1}, $\dot \varepsilon = - \lambda \varepsilon$ during flows, and $\varepsilon^+ = 0$ at jumps.
\end{proof}

Now we have all the ingredients to prove Theorem~\ref{thm:FTClosed_Loop_Stab_HS1}.

\begin{proof}[\unskip\nopunct]{\it Proof of Theorem \ref{thm:FTClosed_Loop_Stab_HS1}:}
Let $\xi = (x,\hat \theta,\psi,\eta,\Phi,Q,m,\tau)$ be a maximal solution to $\HS$ in \eqref{eqn:FTH1} from $\xi(0,0) \in {\cal X}_0$ and consider the following function:
\begin{align} \label{eqn:FTVdef}
    V(t,j) &\coleq \frac{1}{2} \tilde \theta(t,j)^\top \tilde \theta(t,j)
    = \frac{1}{2} |\xi(t,j)|_{\cal A}^2
\end{align}
for all $(t,j) \in \dom \xi$, where $\tilde \theta = \theta - \hat \theta$ is the parameter estimation error.  
Since $\xi(0,0) \in {\cal X}_0$, it follows from Lemma~\ref{prop:FTepsStab} that $\varepsilon(t,j) = 0$ for all $(t,j) \in \dom \xi$. Hence, we have from \eqref{eqn:FTtildethetaC} that
\begin{align} \label{eqn:FTVflow}
    \frac{d}{dt} V(t,0)
    &= - \gamma \tilde \theta(t,0)^\top \psi(t,0)^\top \psi(t,0) \tilde \theta(t,0) \notag \\
    &= - \gamma |\psi(t,0) \tilde \theta(t,0)|^2
    \leq 0
\end{align}
for almost all $t \in I^0 \coleq \{t : (t,0) \in \dom \xi\}$. Thus, $V$ is nonincreasing during the first interval of flow in $\xi$.

Next, we show that $\xi$ converges in finite time to ${\cal A}$.
Since $\xi(0,0) \in {\cal X}_0$, it follows that $\Phi(0,0) = I$ and, thus, \eqref{phi12statetransition} holds for all $t \in I^0$.
Moreover, since $\psi(0,0) = 0$, the solution component $t \mapsto \psi(t,0)$ is the solution to \eqref{eqn:FTpsidynamics0} with $t' = 0$ on $I^0$. Using the fact that $Q(0,0) = 0$, it follows from \eqref{eqn:FTH1} that
\begin{equation*}
    Q(t,0) 
    = \int_{0}^{t} \psi(s,0)^\top \psi(s,0) ds
\end{equation*}
for all $t \in I^0$.
Since Assumptions~\ref{asm:FTphiBound} and~\ref{def:FTFE} hold with $t' = 0$, it follows from Lemma~\ref{lem:FT_FEPsi} that $t \mapsto \psi(t,0)$ is exciting over the interval $[0,T]$, namely, \eqref{eqn:FT_psiFE} holds, and thus $\mu$ in \eqref{eqn:FT_mubar} is positive.
Then, from \eqref{eqn:FT_mubar} with $t' = 0$ and since $\alpha < \mu$, it follows that there exists $t_1 \in (0, T]$ such that $(t_1,0) \in \dom \xi$ and
\begin{equation} \label{QboundProof}
    Q(t_1,0) 
    = \int_{0}^{t_1} \psi(s,0)^\top \psi(s,0) ds
    \geq \alpha I
\end{equation}
and thus $\xi(t_1,0) \in D$. 
Since $\xi$ is maximal, the solution jumps at hybrid time $(t_1,0)$ and, after the jump according to the jump map $G$ in \eqref{eqn:FTH1}, we have from \eqref{hatz_jump} that
\begin{equation} \label{eqn:FThatThetaAtT1}
    \hat \theta(t_1,1) = \theta.
\end{equation}
From the definition of ${\cal A}$ in \eqref{eqn:FTCalA}, we have $\xi(t_1,1)\in {\cal A}$.
Furthermore, since $\xi(t_1,0) \in D$ and $\xi(t_1,1)= G(\xi(t_1,0))$, it follows that $\xi(t_1,1) \in G(D)$, and thus $\xi(t_1,1) \in {\cal A} \cap G(D)$.

We now show that each solution to $\HS$ from ${\cal A} \cap G(D)$ cannot leave ${\cal A}$ via flows or jumps.
Let $\xi^* = (x^*,\hat \theta^*,\psi^*,\eta^*,\Phi^*,Q^*,m^*,\tau^*)$ be a maximal solution to $\HS$ from $\xi^*(0,0) \in {\cal A} \cap G(D)$ and recall the function $V$ in \eqref{eqn:FTVdef}.
%
%
From the definition of $G$ in \eqref{eqn:FTH1}, each $\xi \in \A \cap G(D)$ satisfies $\hat \theta = m$, $\psi = 0$, $\eta = -x$, $\Phi = I$, and $Q = 0$. Hence, $\varepsilon^*(0,0) = 0$ and, by the same arguments as in the proof of Lemma~\ref{prop:FTepsStab}, $\varepsilon^*(t,j) = x^*(t,j) + \eta^*(t,j) - \psi^*(t,j) \theta = 0$ for all $(t,j) \in \dom \xi^*$.
Hence, we have from \eqref{eqn:FTtildethetaC} that, for all $j \in \nats$ and almost all $t \in I^j \coleq \{ t : (t,j) \in \dom \xi^* \}$,
\begin{align} \label{eqn:FTVflow2}
    \frac{d}{dt} V(t,j)
    &= - \gamma \tilde \theta^*(t,j)^{\top} \psi^*(t,j)^{\top} \psi^*(t,j) \tilde \theta^*(t,j) \notag \\
    &= - \gamma |\psi^*(t,j) \tilde \theta^*(t,j)|^2
    \leq 0.
\end{align}
Thus, $V$ is nonincreasing during flows, and $\xi^*$ cannot leave ${\cal A}$ via flow.
We now analyze the variation of $V$ at jumps. Let $(t,j) \in \dom \xi^*$ be such that $(t,j+1) \in \dom \xi^*$.
Since, at the beginning of each interval of flow in $\xi^*$, $\hat \theta^* = m^*$, $\Phi^* = I$, and $\varepsilon^* = 0$, the derivation in \eqref{phi12statetransition}-\eqref{hatz_jump} applies on each interval of flow. Hence, it follows from the jump map $G$ in \eqref{eqn:FTH1} that $\hat \theta^*(t, j+1) = \theta$, and we have
\begin{align} \label{eqn:FTVjump}
    V(t,j+1) - V(t,j) = - V(t,j) \leq 0.
\end{align}
Thus, $V$ is nonincreasing at jumps, and $\xi^*$ cannot leave ${\cal A}$ via a jump. 

In view of \eqref{eqn:FThatThetaAtT1}, \eqref{eqn:FTVflow2}, and \eqref{eqn:FTVjump}, we conclude that, under Assumption~\ref{def:FTFE}, each maximal solution $\xi$ to $\HS$ from $\xi(0,0) \in {\cal X}_0$ converges in finite time to ${\cal A} \cap G(D)$ and remains in ${\cal A}$ for all future hybrid time in $\dom \xi$. Hence, \eqref{eqn:FTFTbound} holds for all $(t,j) \in \dom \xi$ satisfying $t \geq T$ and $j \geq 1$. 
Moreover, by integration using the bounds in \eqref{eqn:FTVflow}, \eqref{eqn:FTVflow2}, and \eqref{eqn:FTVjump}, and from the definition of $V$ in \eqref{eqn:FTVdef}, for each solution $\xi$ to $\HS$ from $\xi(0,0) \in {\cal X}_0$,
\begin{equation*} 
    \frac{1}{2} |\xi(t,j)|_{\cal A}^2 = V(t,j) \leq V(0, 0) = \frac{1}{2} |\xi(0,0)|_{\cal A}^2
\end{equation*}
for all $(t,j) \in \dom \xi$. Hence, \eqref{eqn:FTxibound} holds for all $(t,j) \in \dom \xi$.

To conclude the proof, we bound the condition number of $\Phi-I$ as in \eqref{eqn:FTdeltabound}.
Since Assumption~\ref{asm:FTphiBound} holds and $Q(t_1, 0) \geq \alpha I$ from \eqref{QboundProof}, the arguments in the proof of Lemma~\ref{lem:FTPhiStability} apply with $t' = 0$ and with $T$ replaced by $t_1$.
Then, since $t_1 \leq T$ and $\sigma$ in \eqref{eqn:FTalphabeta} is nondecreasing in $T$, we obtain
\begin{equation} \label{eqn:FTcondOmegaProof}
    \text{\rm cond}(\Omega(t_1,0) - I) \leq \frac{1 + \sigma(\gamma,\alpha)}{1 - \sigma(\gamma,\alpha)}
\end{equation}
with $\sigma$ as in \eqref{eqn:FTalphabeta}, where $\Omega$ is the state transition matrix for \eqref{eqn:FTlemsysa}.
Moreover, since $\Phi(0,0) = I$ and the dynamics during flows of the state component $\Phi$ are equivalent to the dynamics of $\Omega$, it follows that $\Phi(t,0) = \Omega(t,0)$ for all $t \in I^0$, and \eqref{eqn:FTdeltabound} follows from \eqref{eqn:FTcondOmegaProof}.
\end{proof}

\section{Robustness Analysis} \label{sec:FTRobustnessAnalysis}
%
In this section, we analyze the robustness of our algorithm with respect to noise on the state measurements.
Consider additive noise $t \mapsto \nu(t) \in \reals^n$ in the measurements of the plant state $x$.
During flows, the measurement noise affects the dynamics of the state components $\hat \theta$ and $\eta$ in \eqref{eqn:FTH1}. In particular, the dynamics of $\hat \theta$ become
\begin{equation*}
    \dot{\hat \theta} = \gamma \psi^\top (x + \nu(t) + \eta - \psi \hat \theta)
\end{equation*}
and the dynamics of $\eta$ become
\begin{equation*}
    \dot{\eta} = -\lambda (x + \nu(t) + \eta) - f(x + \nu(t), u(t)).
\end{equation*}
Similarly, at jumps the measurement noise affects the jump map for $\eta$, which becomes
\begin{equation*}
    \eta^+ = -x - \nu(t)
\end{equation*}
which implies that $\varepsilon$ in \eqref{eqn:FTeps} is reset as $\varepsilon^+ = -\nu(t)$ when a jump occurs at hybrid time $(t,j)$.

We denote the hybrid system $\HS$ in \eqref{eqn:FTH1} under the effect of the measurement noise $\nu$ as $\HS_\nu$, with state $\xi \coleq (x, \hat \theta, \psi, \eta, \Phi, Q, m, \tau) \in {\cal X}$ and dynamics
\begin{equation} \label{eqn:FTHSnu}
\begin{aligned} 
    {\cal H}_\nu : 
    \left\{
    \begin{aligned}
    \dot \xi &= F_\nu(\xi) & \qquad
    \xi &\in C_\nu \\
    \xi^+ &= G_\nu(\xi) & \qquad
    \xi &\in D_\nu
    \end{aligned}
    \right.
\end{aligned} 
\end{equation}
where
\begin{equation} \label{eqn:FTFnuGnu}
\mbox{$\squeezespaces{0.40}
\begin{aligned}
    F_\nu(\xi) &\coleq 
    \begin{bmatrix}
        f(x,u(\tau)) + \phi(\tau) \theta \\
        \gamma \psi^\top (x + \nu(\tau) + \eta - \psi \hat \theta) \\
        -\lambda \psi + \phi(\tau) \\
        -\lambda (x + \nu(\tau) + \eta) - f(x + \nu(\tau), u(\tau)) \\
        - \gamma \psi^\top \psi \Phi \\
        \psi^\top \psi \\
        0 \\
        1
    \end{bmatrix} &
    &\xi \in C_\nu \\
    G_\nu(\xi) &\coleq
    (x, R(\xi), 0, -x-\nu(\tau), I, 0, R(\xi), \tau) &
    &\xi \in D_\nu
\end{aligned}
$}
\end{equation}
with $R$ as in \eqref{eqn:FTR} and $C_\nu \coleq C$, $D_\nu \coleq D$, with $C$, $D$ as in \eqref{eqn:FTCandD}.

To enable analysis of the noise effect, we make the following Lipschitz continuity assumption.
\begin{assumption} \label{asm:FTboundedf}
%
%
Given the system \eqref{eqn:FTxdot}, there exists $L_x > 0$ such that, for all $x_1, x_2 \in \reals^n$ and all $u \in \reals^r$,
\begin{equation*}
    |f(x_1,u) - f(x_2,u)| \leq L_x |x_1 - x_2|.
\end{equation*}
\end{assumption}
\subsection{Robustness Under Persistence of Excitation}
It can be shown that the conditions of Theorem~\ref{thm:FTClosed_Loop_Stab_HS1} -- namely, excitation over only a finite interval as in Assumption~\ref{def:FTFE} -- ensure that each maximal solution to the hybrid system $\HS_\nu$ from an appropriate set of initial conditions jumps at least once. However, unlike in Theorem~\ref{thm:FTClosed_Loop_Stab_HS1}, the parameter estimation error for $\HS_\nu$ is not guaranteed to converge to zero after the first jump due to the noise $\nu$ in the measurements of $x$.
Hence, we impose the following PE condition to ensure that each maximal solution to $\HS_\nu$ jumps periodically, which allows us to establish an ISS bound on the estimation error in Theorem~\ref{thm:FTISS} below.
\begin{assumption} \label{def:FTPE}
Given $t \mapsto \phi(t) \in \reals^{n \times p}$, there exist $T, \mu_\phi > 0$, and $L_\phi \geq 0$ such that
\begin{equation}  \label{eqn:FTPE}
    \int_{t'}^{t'+T} 
    \phi (s)^\top 
    \phi (s) ds
    \geq \mu_\phi I
    \qquad \forall t' \geq 0
\end{equation}
%
%
and, for all $s, v \geq 0$,
\begin{equation}  \label{eqn:LischitzPE_phi}
    |\phi(s) - \phi(v)| \leq L_\phi |s - v|.
\end{equation}
\end{assumption}

%
\begin{remark}
Compared to Assumption~\ref{def:FTFE}, which requires that the inequalities in \eqref{eqn:FTFE} and \eqref{eqn:Lischitz_phi} hold over only a finite interval, Assumption~\ref{def:FTPE} requires that \eqref{eqn:FTPE} holds for all $t' \geq 0$ and that \eqref{eqn:LischitzPE_phi} holds for all $s, v \geq 0$.
%
\end{remark}
%
%
\begin{theorem} \label{thm:FTISS}
Given $t \mapsto \phi(t) \in \reals^{n \times p}$, $t \mapsto u(t) \in \reals^r$, $(x,u) \mapsto f(x,u) \in \reals^n$ in \eqref{eqn:FTxdot}, $t \mapsto \nu(t) \in \reals^{n}$, 
and $\lambda, \alpha > 0$, suppose that Assumptions~\ref{asm:FTphiBound},~\ref{asm:FTboundedf}, and~\ref{def:FTPE} hold with $T, \phi_M, \mu_\phi > 0$, $L_\phi \geq 0$, and $\alpha < \mu$ for all $t' \geq 0$, with $\mu$ as in \eqref{eqn:FT_mubar}. 
Then, for each $\gamma > 0$ and each $\zeta \in (0, 1)$, the hybrid system $\HS_\nu$ in \eqref{eqn:FTHSnu} is such that the parameter estimation error $(t,j) \mapsto \tilde \theta(t,j) \coleq \theta - \hat \theta(t,j)$ for each solution $\xi$ to $\HS_\nu$ from $\xi(0,0) \in {\cal X}_{\nu,0} \coleq \{\xi \in {\cal X} : \hat \theta = m, \, \psi = 0, \, \eta = -x-\nu(0), \, \Phi = I, \, Q = 0, \, \tau = 0 \}$ satisfies
%
%
%
\begin{equation} \label{eqn:FTxiboundISS}
\begin{aligned}
    |\tilde \theta(t,j)|
    &\leq q(j) |\tilde \theta(0,0)|
    + (1 - q(j)) \kappa_1(\|\nu\|_\infty) \\
    &\quad + \kappa_2 (\|\nu\|_\infty)
\end{aligned}
\end{equation}
for all $(t,j) \in \dom \xi$, where  
$q(j) \coleq 1$ {\rm if} $j = 0$, $q(j) \coleq 0$ {\rm if} $j \in \nats \setminus \{0\}$, and
\begin{subequations} \label{eqn:FTGammaOmegaRho}
\begin{align} 
    \kappa_1(s) &\coleq \frac{\gamma \psi_M}{1 - \sigma(\gamma,\alpha)} \left( \frac{1}{\omega} + T (\rho + 1) \right) s \label{eqn:FTalpha1} \\
    \kappa_2(s) &\coleq \gamma \psi_M \left( \frac{1}{\omega} + T (\rho + 1) \right) s \label{eqn:FTalpha2} \\
    \omega &\coleq \frac{1}{2} \min \{ \lambda (1 - \zeta), -\ln(\zeta) \} \label{eqn:FTomega} \\
    \rho &\coleq \max \left\{ \frac{2 (L_x + \lambda)}{\lambda \sqrt{\zeta}} , \frac{1}{\sqrt{\zeta}} \right\}, \label{eqn:FTrho}
\end{align}
\end{subequations}
with $\sigma, \psi_M$ as in \eqref{eqn:FTalphabeta} and $L_x$ from Assumption~\ref{asm:FTboundedf}.
\end{theorem}
\startmodifnew
%
\stopmodifnew

The proof of Theorem~\ref{thm:FTISS} proceeds as follows. 
First, in Lemma~\ref{lem:FTHSnuSol}, we establish expressions for the parameter estimation error during each interval of flow and immediately after each jump for each solution to $\HS_{\nu}$ from ${\cal X}_{\nu,0}$.
Next, in Lemma~\ref{lem:FTpsiHSbound}, we establish a uniform upper bound on the norm of the state component $\psi$ for each solution to $\HS_{\nu}$ from ${\cal X}_{\nu,0}$. 
Then, in Lemma~\ref{lem:FTepsISS}, we establish an upper bound on the norm of the quantity $\varepsilon$ in \eqref{eqn:FTeps} for each solution to $\HS_{\nu}$. 
Next, in Lemma~\ref{lem:FTPEtoFE}, we establish that the PE condition in Assumption~\ref{def:FTPE} implies the FE condition in Assumption~\ref{def:FTFE}.
Finally, these results are used to prove Theorem~\ref{thm:FTISS}.
\begin{lemma} \label{lem:FTHSnuSol}
Given $t \mapsto \phi(t) \in \reals^{n \times p}$, $t \mapsto u(t) \in \reals^r$, $(x,u) \mapsto f(x,u) \in \reals^n$ in \eqref{eqn:FTxdot}, $t \mapsto \nu(t) \in \reals^{n}$, 
the hybrid system $\HS_\nu$ in \eqref{eqn:FTHSnu} is such that, for each solution $\xi$ to $\HS_\nu$ from $\xi(0,0) \in {\cal X}_{\nu,0}$, the parameter estimation error $(t,j) \mapsto \tilde \theta(t,j) \coleq \theta - \hat \theta(t,j)$ satisfies
\begin{align} \label{eqn:FTTildeThetaFlow}
    \tilde \theta(t,j) 
    &= \Phi(t,j) \tilde \theta(t_j,j) \\
    &\quad - \gamma \int_{t_j}^{t} \Phi(s,j) \psi(s,j)^\top( \varepsilon(s,j) + \nu(s)) ds \notag
\end{align}
for each $j \in \nats$ and all $t \in I^j \coleq \{t : (t,j) \in \dom \xi\}$, where $\{t_j\}_{j = 0}^{\sup_j \dom \xi}$ is the sequence defining $\dom \xi$ as in Section~\ref{sec:COMHybridSystems}, and $\varepsilon$ is defined in \eqref{eqn:FTeps}. 
Moreover, for each $(t,j) \in \dom \xi$ such that $(t,j+1) \in \dom \xi$,
\begin{align} \label{eqn:FTTildeThetaJump}
        &\tilde \theta(t,j+1) \\
        &\quad = - K_1(\xi(t,j)) \gamma \int_{t_j}^{t} \Phi(s,j) \psi(s,j)^\top( \varepsilon(s,j) + \nu(s)) ds \notag
\end{align} 
with $K_1$ as in \eqref{eqn:FTKi1}.
\end{lemma}
\begin{proof}
\iftoggle{report}{
    This proof is in Appendix~\ref{apx:FTHSnuSol}.
}{
    This proof is in \cite{JohnsonFiniteTimeArxiv}.
}
\end{proof}
\begin{lemma} \label{lem:FTpsiHSbound}
Given $t \mapsto \phi(t) \in \reals^{n \times p}$, $t \mapsto u(t) \in \reals^r$, $(x,u) \mapsto f(x,u) \in \reals^n$ in \eqref{eqn:FTxdot}, $t \mapsto \nu(t) \in \reals^{n}$, 
suppose that Assumption~\ref{asm:FTphiBound} holds. Then, for each $\lambda > 0$, the hybrid system $\HS_\nu$ in \eqref{eqn:FTHSnu} is such that the $\psi$ component of each solution $\xi$ to $\HS_\nu$ from $\xi(0,0) \in {\cal X}_{\nu,0}$ satisfies
\begin{equation} \label{eqn:FTpsiHSbound}
    |\psi(t,j)| \leq \psi_M
\end{equation}
for all $(t,j) \in \dom \xi$, with $\psi_M$ as in \eqref{eqn:FTalphabeta}.
\end{lemma}
\begin{proof} 
The proof of Lemma~\ref{lem:FTpsiHSbound} follows directly from Lemma~\ref{lem:FTPhiExcitation} and the fact that $\psi(0,0) = 0$ and, at each jump, $\psi^+ = 0$.
\end{proof}
\begin{lemma} \label{lem:FTepsISS}
Given $t \mapsto \phi(t) \in \reals^{n \times p}$, $t \mapsto u(t) \in \reals^r$, $(x,u) \mapsto f(x,u) \in \reals^n$ in \eqref{eqn:FTxdot}, $t \mapsto \nu(t) \in \reals^{n}$, 
suppose that Assumption~\ref{asm:FTboundedf} holds. Then, for each $\lambda,\gamma > 0$ and each $\zeta \in (0,1)$, the hybrid system $\HS_\nu$ in \eqref{eqn:FTHSnu} is such that, for each solution $\xi$ to $\HS_\nu$, $(t,j)~\mapsto~\varepsilon(t,j) = x(t,j) + \eta(t,j) - \psi(t,j) \theta$ in \eqref{eqn:FTeps} satisfies
\begin{equation} \label{eqn:FTepsISS}
    |\varepsilon(t,j)| \leq \mathrm{e}^{-\omega (t+j-t'-j')} |\varepsilon(t',j')| + \rho \|\nu\|_\infty
\end{equation}
for all $(t,j),(t',j') \in \dom \xi$ satisfying $t+j \geq t'+j'$, with $\omega$ as in \eqref{eqn:FTomega} and $\rho$ as in \eqref{eqn:FTrho}.
\end{lemma}
\iftoggle{report}{
\begin{proof} 
    This proof is in Appendix~\ref{apx:FTepsISS}.
\end{proof} 
}{
\begin{proof} 
    This proof is in \cite{JohnsonFiniteTimeArxiv}.
\end{proof} 
}
\begin{lemma} \label{lem:FTPEtoFE}
Given $t \mapsto \phi(t) \in \reals^{n \times p}$, suppose that Assumption~\ref{def:FTPE} holds. Then, for all $t' \geq 0$, Assumption~\ref{def:FTFE} holds with $T, \mu_\phi, L_\phi$ from Assumption~\ref{def:FTPE}.
\end{lemma}
\begin{proof} 
    The proof of Lemma~\ref{lem:FTPEtoFE} is immediate since, by Assumption~\ref{def:FTPE}, \eqref{eqn:FTFE} and \eqref{eqn:Lischitz_phi} hold for all $t' \geq 0$.
\end{proof}
Now we have all the ingredients to prove Theorem~\ref{thm:FTISS}.
\begin{proof}[\unskip\nopunct]{\it Proof of Theorem~\ref{thm:FTISS}:}
Let $\xi = (x, \hat \theta, \psi, \eta, \Phi, Q, m, \tau)$ be a maximal solution to $\HS_\nu$ in \eqref{eqn:FTHSnu} from $\xi(0,0) \in {\cal X}_{\nu,0}$.
We first bound the parameter estimation error $\tilde \theta = \theta - \hat \theta$ during the initial interval of flow $I^0 \coleq \{t : (t,0) \in \dom \xi\}$. 
From Lemma~\ref{lem:FTHSnuSol}, we have that for all $t \in I^0$,
\begin{align*} 
    \tilde \theta(t,0)
    &= \Phi(t,0) \tilde \theta(0,0) \\
    &\quad- \gamma \int_{0}^{t} \Phi(s,0) \psi(s,0)^\top( \varepsilon(s,0) + \nu(s)) ds.
\end{align*}
Since Assumption~\ref{def:FTPE} is satisfied, it follows from Lemma~\ref{lem:FTPEtoFE} that Assumption~\ref{def:FTFE} holds over the interval $[0, T]$ and, since Assumption~\ref{asm:FTphiBound} holds, the conditions of Lemma~\ref{lem:FTPhiStability} are satisfied with $t' = 0$. 
Moreover, since $\Phi(0,0) = I$ and the dynamics during flows of the state component $\Phi$ are equivalent to the dynamics of the state transition matrix $\Omega$ for \eqref{eqn:FTlemsysa}, it follows from \eqref{eqn:FTOmegaUpperBound} that $|\Phi(t,0)| \leq 1$ for all $t \in I^0$. Hence,
\begin{align} \label{eqn:FTISSErrorBoundI0}
    &|\tilde \theta(t,0)|
    \leq | \Phi(t,0) | | \tilde \theta(0,0) | \notag \\
    &\quad + \gamma \int_{0}^{t} |\Phi(s,0)| |\psi(s,0)| ( |\varepsilon(s,0)| + |\nu(s)| ) ds \\
    &\leq | \tilde \theta(0,0) | + \gamma \psi_M \int_{0}^{t} \big( \mathrm{e}^{-\omega s} |\varepsilon(0,0)| + (\rho + 1) \|\nu\|_\infty \big) ds \notag \\
    %
    %
    &\leq |\tilde \theta(0,0)| 
    + \gamma \psi_M \bigg( \frac{1}{\omega} |\varepsilon(0,0)|
    + t (\rho + 1) \|\nu\|_\infty \bigg) \notag \\
    &\leq |\tilde \theta(0,0)| 
    + \gamma \psi_M \left( \frac{1}{\omega} + t(\rho + 1) \right) \|\nu\|_\infty \notag
\end{align}
for all $t \in I^0$,
which follows from Lemmas~\ref{lem:FTpsiHSbound} and~\ref{lem:FTepsISS}, 
the fact that $\omega > 0$, 
and the fact that, since $\xi(0,0) \in {\cal X}_{\nu,0}$, $|\varepsilon(0,0)| = |\nu(0)| \leq \|\nu\|_\infty$.

Note that the measurement noise $\nu$ does not affect the dynamics of the $\psi$ or $Q$ components of the solution $\xi$ of $\HS_\nu$. Hence, by Assumptions~\ref{asm:FTphiBound} and~\ref{def:FTPE} and Lemma~\ref{lem:FTPEtoFE}, the conditions of Theorem~\ref{thm:FTClosed_Loop_Stab_HS1} hold. Moreover, by the equivalence between the jump sets of $\HS$ in \eqref{eqn:FTH1} and $\HS_\nu$ in \eqref{eqn:FTHSnu}, it follows from Theorem~\ref{thm:FTClosed_Loop_Stab_HS1} that $\xi$ jumps at hybrid time $(t_1,0)$ with $t_1 \in (0,T]$, and thus $I^0 \subset [0, T]$. Hence, for all $t \in I^0$,
\begin{align} \label{eqn:FTHSnuBound1}
    |\tilde \theta(t,0)| 
    %
    &\leq |\tilde \theta(0,0)| + \kappa_2(\|\nu\|_\infty)
\end{align}
with $\kappa_2$ as in \eqref{eqn:FTalpha2}.

Next, we bound the parameter estimation error after each jump in $\dom \xi$. Let $(t_j,j) \in \dom \xi$ be such that $(t_j,j-1) \in \dom \xi$. From Lemma~\ref{lem:FTHSnuSol}, we have
\begin{equation*}
    \begin{aligned}
        &\tilde \theta(t_j,j) = - K_1(\xi(t_j,j-1)) \gamma \\
        &\quad \times \int_{t_{j-1}}^{t_j} \Phi(s,j-1) \psi(s,j-1)^\top( \varepsilon(s,j-1) + \nu(s)) ds
    \end{aligned}
\end{equation*} 
Since Assumption~\ref{def:FTPE} is satisfied, it follows from Lemma~\ref{lem:FTPEtoFE} that Assumption~\ref{def:FTFE} holds over the interval $[t_{j-1}, t_{j-1}+T]$ and, since Assumption~\ref{asm:FTphiBound} holds, the conditions of Lemma~\ref{lem:FTPhiStability} are satisfied with $t' = t_{j-1}$.
Moreover, since $\Phi(t_{j-1},j-1) = I$ and the dynamics during flows of $\Phi$ in \eqref{eqn:FTHSnu} are equivalent to the dynamics of the state transition matrix $\Omega$ for \eqref{eqn:FTlemsysa}, it follows from \eqref{eqn:FTOmegaUpperBound} that $|\Phi(t,j-1)| \leq 1$ for all $t \in I^{j-1}$. Thus,
\begin{align} \label{eqn:FTHsnuJumpBound}
        &|\tilde \theta(t_j,j)|
        \leq |K_1(\xi(t_j,j-1))| \gamma \\
        &\quad \mbox{$\squeezespaces{0.50} \times \displaystyle\int_{t_{j-1}}^{t_j} |\Phi(s,j-1)| |\psi(s,j-1)| ( |\varepsilon(s,j-1)| + |\nu(s)| ) ds $} \notag \\
        &\leq |K_1(\xi(t_j,j-1))| \gamma \psi_M \notag \\
        &\quad \times \int_{t_{j-1}}^{t_j} \big( \mathrm{e}^{-\omega (s+j-1)} |\varepsilon(0,0)| + (\rho + 1) \|\nu\|_\infty \big) ds \notag \\
        %
        %
        &\leq |K_1(\xi(t_j,j-1))| \gamma \psi_M \notag \\
        &\quad \times \bigg( \frac{1}{\omega} |\varepsilon(0,0)|
        + (t_j-t_{j-1}) (\rho + 1) \|\nu\|_\infty \bigg) \notag \\
        &\leq \mbox{$\squeezespaces{0.7} |K_1(\xi(t_j,j-1))| \gamma \psi_M
        \left( \displaystyle\frac{1}{\omega} + (t_j-t_{j-1}) (\rho + 1) \right) \|\nu\|_\infty $} \notag
\end{align} 
which follows from Lemmas~\ref{lem:FTpsiHSbound} and~\ref{lem:FTepsISS},
the fact that $\omega > 0$,
and the fact that $|\varepsilon(0,0)| = |\nu(0)| \leq \|\nu\|_\infty$.
Since, by Assumption~\ref{def:FTPE} and Lemma~\ref{lem:FTPEtoFE}, Assumption~\ref{def:FTFE} holds over the interval $[t_{j-1}, t_{j-1}+T]$, it follows from the same arguments as in the proof of Theorem~\ref{thm:FTClosed_Loop_Stab_HS1} that $t_j - t_{j-1} \leq T$.
Moreover, since $Q(t_j, j-1) \geq \alpha I$, the arguments in the proof of Lemma~\ref{lem:FTPhiStability} apply with $t' = t_{j-1}$ and with $T$ replaced by $t_j - t_{j-1}$ and, since $\sigma$ in \eqref{eqn:FTalphabeta} is nondecreasing in $T$, \eqref{eqn:FTK1Bound} holds for $\Phi(t_j, j-1)$.
Hence,
\begin{align} \label{eqn:FTHSnuBound01}
        &|\tilde \theta(t_j,j)| \\
        &\leq |K_1(\xi(t_j,j-1))| \gamma \psi_M \left( \frac{1}{\omega} + T (\rho + 1) \right) \|\nu\|_\infty \notag \\
        &\leq \frac{\gamma \psi_M}{1 - \sigma(\gamma,\alpha)} \left( \frac{1}{\omega} + T (\rho + 1) \right) \|\nu\|_\infty 
        = \kappa_1(\|\nu\|_\infty) \notag
\end{align} 
with $\kappa_1$ as in \eqref{eqn:FTalpha1} and $\sigma$ as in \eqref{eqn:FTalphabeta}, where the last inequality follows from the definition of $K_1 = - ( \Phi - I )^{-1}$ in \eqref{eqn:FTKi1} and from \eqref{eqn:FTK1Bound}.

During the interval of flow after the jump at hybrid time $(t_j,j)$, we have from Lemma~\ref{lem:FTHSnuSol} that, for all $t \in I^j$,
\begin{align*} 
    \tilde \theta(t,j) 
    &= \Phi(t,j) \tilde \theta(t_j,j) \\
    &\quad - \gamma \int_{t_j}^{t} \Phi(s,j) \psi(s,j)^\top( \varepsilon(s,j) + \nu(s)) ds.
\end{align*}
Since Assumption~\ref{def:FTPE} is satisfied, it follows from Lemma~\ref{lem:FTPEtoFE} that Assumption~\ref{def:FTFE} holds over the interval $[t_j, t_j+T]$ and, since Assumption~\ref{asm:FTphiBound} holds, the conditions of Lemma~\ref{lem:FTPhiStability} are satisfied with $t' = t_j$. 
Moreover, since $\Phi(t_{j},j) = I$ and the dynamics during flows of $\Phi$ in \eqref{eqn:FTHSnu} are equivalent to the dynamics of the state transition matrix $\Omega$ for \eqref{eqn:FTlemsysa}, it follows from \eqref{eqn:FTOmegaUpperBound} that $|\Phi(t,j)| \leq 1$ for all $t \in I^{j}$. 
Using similar arguments as in \eqref{eqn:FTISSErrorBoundI0}, we obtain
\begin{align*}
    |\tilde \theta(t,j)| 
    %
    %
    %
    %
    %
    %
    &\leq |\tilde \theta(t_j,j)| 
    + \gamma \psi_M \left( \frac{1}{\omega} + (t - t_j)(\rho + 1) \right) \|\nu\|_\infty 
\end{align*}
for all $t \in I^j$. 
Since, by Assumption~\ref{def:FTPE} and Lemma~\ref{lem:FTPEtoFE}, Assumption~\ref{def:FTFE} holds over the interval $[t_{j}, t_{j}+T]$, it follows from the same arguments as in the proof of Theorem~\ref{thm:FTClosed_Loop_Stab_HS1} that $t_{j+1} - t_{j} \leq T$. Hence, for all $t \in I^j$,
\begin{equation} \label{eqn:FTHSnuBound2}
\begin{aligned} 
    |\tilde \theta(t,j)| 
    &\leq |\tilde \theta(t_j,j)| 
    + \gamma \psi_M \left( \frac{1}{\omega} + T(\rho + 1) \right) \|\nu\|_\infty \\
    &\leq \kappa_1(\|\nu\|_\infty) + \kappa_2(\|\nu\|_\infty)
\end{aligned}
\end{equation}
where the last line follows from substituting \eqref{eqn:FTHSnuBound01}.

Since the bound in \eqref{eqn:FTHSnuBound2} does not depend on the choice of hybrid time instant $(t_j,j)$, it follows that \eqref{eqn:FTHSnuBound2} holds for each $j \in \nats \setminus \{0\}$ and all $t \in I^j$. 
Hence, we unify the bounds in \eqref{eqn:FTHSnuBound1} and \eqref{eqn:FTHSnuBound2} to obtain that
\begin{equation*}
    \mbox{$\squeezespaces{0.6}
    |\tilde \theta(t,j)| 
    \leq q(j) |\tilde \theta(0,0)|
    + (1 - q(j)) \kappa_1(\|\nu\|_\infty)
    + \kappa_2(\|\nu\|_\infty)
    $}
\end{equation*}
for all $(t,j) \in \dom \xi$, where $q(j) \coleq 1$ if $j = 0$ and $q(j) \coleq 0$ if $j \in \nats \setminus \{0\}$. Hence, \eqref{eqn:FTxiboundISS} holds.
%
{}
\end{proof}

\subsection{Robustness Under Finite Excitation}
To analyze the robustness properties of $\HS_\nu$ under the finite excitation condition in Assumption~\ref{def:FTFE}, we modify the algorithm to jump at most once by augmenting the state vector with a new component $k \in \{0, 1\}$. The resulting algorithm, denoted by ${\HS}_\nu',$ has state $\xi' \coleq (x, \hat \theta, \psi, \eta, \Phi, Q, m, \tau, k) \in {\cal X}' \coleq {\cal X} \times \{0, 1\}$ and dynamics
\begin{equation} \label{eqn:FTHSnu2}
\begin{aligned} 
    {\HS}_\nu' : 
    \left\{
    \begin{aligned}
    \dot \xi' &=
    \begin{bmatrix}
        F_\nu(\xi') \\
        0
    \end{bmatrix} =:
    F_\nu'(\xi') & \qquad
    \xi' &\in C_\nu' \\
    \xi'^+ &=
    \begin{bmatrix}
        G_\nu(\xi') \\
        1
    \end{bmatrix} =:
    G_\nu'(\xi') & \qquad
    \xi' &\in D_\nu'
    \end{aligned}
    \right.
\end{aligned} 
\end{equation}
with $F_\nu, G_\nu$ as in \eqref{eqn:FTFnuGnu},
\begin{align} \label{eqn:FTCandD2}
    C_\nu' &\coleq \left\{ \xi' \in {\cal X}' : \lambda_{\min}(Q) \leq \alpha \right\} \cup
    \left\{ \xi' \in {\cal X}' : k = 1 \right\} \notag \\
    D_\nu' &\coleq \left\{ \xi' \in {\cal X}' : \lambda_{\min}(Q) \geq \alpha, \, k = 0 \right\},
\end{align}
and $\alpha > 0$ is a design parameter. 
The jump set $D_\nu'$ in \eqref{eqn:FTCandD2} ensures that ${\HS}_\nu'$ jumps only once, when $k = 0$. 
Note that the results in Lemmas~\ref{lem:FTHSnuSol},~\ref{lem:FTpsiHSbound}, and~\ref{lem:FTepsISS} also apply to the hybrid system ${\HS}_\nu'$.

We now show that, under the finite excitation condition in Assumption~\ref{def:FTFE}, the hybrid system ${\HS}_\nu'$ is such that, for solutions from an augmented version of ${\cal X}_{\nu,0}$ defined above \eqref{eqn:FTxiboundISS}, the norm of the parameter estimation error is upper bounded by a function of the integral of the noise $\nu$ in the measurements of $x$.
\begin{theorem} \label{thm:FTiISS}
Given $t \mapsto \phi(t) \in \reals^{n \times p}$, $t \mapsto u(t) \in \reals^r$, $(x,u) \mapsto f(x,u) \in \reals^n$ in \eqref{eqn:FTxdot}, $t \mapsto \nu(t) \in \reals^{n}$, 
and $\lambda, \alpha > 0$, suppose that Assumptions~\ref{asm:FTphiBound},~\ref{def:FTFE}, and~\ref{asm:FTboundedf} hold with $t' = 0$, $T, \phi_M, \mu_\phi > 0$, $L_\phi \geq 0$, and $\alpha < \mu$ with $\mu$ as in \eqref{eqn:FT_mubar}. 
Then, for each $\gamma > 0$ and each $\zeta \in (0, 1)$, the hybrid system ${\HS}_\nu'$ in \eqref{eqn:FTHSnu2} is such that the parameter estimation error $(t,j) \mapsto \tilde \theta(t,j) \coleq \theta - \hat \theta(t,j)$ for each solution $\xi'$ to ${\HS}_\nu'$ from $\xi'(0,0) \in {\cal X}_{\nu,0}' \coleq {\cal X}_{\nu,0} \times \{0\}$ satisfies
\begin{align} \label{eqn:FTxiboundiISS}
    |\tilde \theta(t,j)|
    &\leq q(j) |\tilde \theta(0,0)|
    + (1 - q(j)) \kappa_1'(\|\nu\|_\infty) \\
    &\qquad + \kappa_2'(\|\nu\|_\infty) 
    + \int_{0}^{t} \kappa_3'(|\nu(s)|,\|\nu\|_\infty) ds \notag
\end{align}
for all $(t,j) \in \dom \xi'$, where 
$q(j) \coleq 1$ {\rm if} $j = 0$, $q(j) \coleq 0$ {\rm if} $j \in \nats \setminus \{0\}$, and 
\begin{align} \label{eqn:FTcalV}
    \kappa_1'(s) &\coleq \frac{\gamma \psi_M}{1 - \sigma(\gamma,\alpha)} \left( \frac{1}{\omega} + T (\rho + 1) \right) s \notag \\
    \kappa_2'(s) &\coleq \frac{\gamma \psi_M}{\omega} s \\
    \kappa_3'(s_1, s_2) &\coleq \gamma \psi_M \left( s_1 + \rho s_2 \right) \notag
\end{align}
%
%
%
with $\psi_M$ as in \eqref{eqn:FTalphabeta} and $\omega, \rho$ as in \eqref{eqn:FTGammaOmegaRho}.
\end{theorem}
The proof of Theorem~\ref{thm:FTiISS} follows along the same lines as the proof of Theorem~\ref{thm:FTISS}. However, since each solution $\xi'$ to $\HS_\nu'$ in \eqref{eqn:FTHSnu2} jumps at most once, we cannot leverage periodic jumps in $\xi'$ to bound the effect of the measurement noise.
\iftoggle{test}{
\begin{proof}
    Let $\xi' = (x, \hat \theta, \psi, \eta, \Phi, Q, m, \tau, k)$ be a maximal solution to ${\HS}_\nu'$ in \eqref{eqn:FTHSnu2} from $\xi'(0,0) \in {\cal X}_{\nu,0}'$.
    We first bound the parameter estimation error $\tilde \theta = \theta - \hat \theta$ during the initial interval of flow $I^0 \coleq \{t : (t,0) \in \dom \xi'\}$. 
    From Lemma~\ref{lem:FTHSnuSol}, we have that for all $t \in I^0$,
    \begin{align*} 
        \tilde \theta(t,0)
        &= \Phi(t,0) \tilde \theta(0,0) - \gamma \int_{0}^{t} \Phi(s,0) \psi(s,0)^\top( \varepsilon(s,0) + \nu(s)) ds.
    \end{align*}
    Since $\Phi(0,0) = I$ and the dynamics during flows of $\Phi$ in \eqref{eqn:FTHSnu2} are equivalent to the dynamics of the state transition matrix $\Omega$ for \eqref{eqn:FTlemsysa}, it follows from \eqref{eqn:FTOmegaUpperBound} that $|\Phi(t,0)| \leq 1$ for all $t \in I^{0}$. Thus,
    \begin{align}
        |\tilde \theta(t,0)| 
        &\leq | \Phi(t,0) | | \tilde \theta(0,0) | + \gamma \int_{0}^{t} |\Phi(s,0)| |\psi(s,0)| ( |\varepsilon(s,0)| + |\nu(s)| ) ds \notag \\
        &\leq | \theta(0,0) | + \gamma \psi_M \int_{0}^{t} \big( \mathrm{e}^{-\omega s} |\varepsilon(0,0)| 
        + \int_{0}^{t} \left( |\nu(s)| + \rho \|\nu\|_\infty \right) ds \bigg) \notag \\
        &\leq |\tilde \theta(0,0)| + \gamma \psi_M \bigg( \int_{0}^{\infty} \mathrm{e}^{- \omega s} |\varepsilon(0,0)| ds
        + \int_{0}^{t} \left( |\nu(s)| + \rho \|\nu\|_\infty \right) ds \bigg) \notag \\
        &\leq |\tilde \theta(0,0)| 
        + \gamma \psi_M \bigg( \frac{1}{\omega} |\varepsilon(0,0)|
        + \int_{0}^{t} \left( |\nu(s)| + \rho \|\nu\|_\infty \right) ds \bigg) \notag \\
        &\leq |\tilde \theta(0,0)| 
        + \gamma \psi_M \bigg( \frac{1}{\omega} \|\nu\|_\infty
        + \int_{0}^{t} \left( |\nu(s)| + \rho \|\nu\|_\infty \right) ds \bigg) \notag \\
        &\leq |\tilde \theta(0,0)|
        + \kappa_2'(\|\nu\|_\infty)
        + \int_{0}^{t} \kappa_3'(|\nu(s)|,\|\nu\|_\infty) ds \label{eqn:FTHSnuBound012}
    \end{align}
    for all $t \in I^0$, with $\kappa_2',\kappa_3'$ as in \eqref{eqn:FTcalV}, where the second line follows from Lemmas~\ref{lem:FTpsiHSbound} and~\ref{lem:FTepsISS},
    the fourth line uses the fact that $\omega > 0$, 
    and the fifth line follows from the fact that $|\varepsilon(0,0)| = |\nu(0)| \leq \|\nu\|_\infty$.
    
    Note that the measurement noise $\nu$ does not affect the dynamics of the $\psi$ or $Q$ components of the solution $\xi'$ of $\HS_\nu'$. Hence, by Assumptions~\ref{asm:FTphiBound} and~\ref{def:FTFE}, the conditions of Theorem~\ref{thm:FTClosed_Loop_Stab_HS1} hold. Thus, it follows from the same arguments as in the proof of Theorem~\ref{thm:FTClosed_Loop_Stab_HS1} that 
    the solution $\xi'$ jumps at hybrid time $(t_1,0)$ with $t_1 \in (0, t'+T]$. After the jump, we have from Lemma~\ref{lem:FTHSnuSol} that 
    \begin{equation*}
        \begin{aligned}
            \tilde \theta(t_1,1) = - K_1(\xi'(t_1,0)) \gamma \int_{0}^{t_1} \Phi(s,0) \psi(s,0)^\top( \varepsilon(s,0) + \nu(s)) ds
        \end{aligned}
    \end{equation*} 
    which is upper bounded by
    \begin{align*}
        |\tilde \theta(t_1,1)|
        &\leq |K_1(\xi'(t_1,0))| \gamma \int_{0}^{t_1} |\Phi(s,0)| |\psi(s,0)| ( |\varepsilon(s,0)| + |\nu(s)| ) ds \notag \\
        &\leq |K_1(\xi'(t_1,0))| \gamma \psi_M \int_{0}^{t_1} \big( \mathrm{e}^{-\omega s} |\varepsilon(0,0)| 
        + \int_{0}^{t_1} \left( 1 + \rho \right) \|\nu\|_\infty ds \bigg) \notag \\
        &\leq |K_1(\xi'(t_j,j-1))| \gamma \psi_M \bigg( \int_{0}^{\infty} \mathrm{e}^{- \omega s} |\varepsilon(0,0)| ds
        + t_1 \left( 1 + \rho \right) \|\nu\|_\infty \bigg) \notag \\
        &\leq |K_1(\xi'(t_j,j-1))| \gamma \psi_M \bigg( \frac{1}{\omega} |\varepsilon(0,0)|
        + t_1 \left( 1 + \rho \right) \|\nu\|_\infty \bigg) \notag \\
        &\leq |K_1(\xi'(t_j,j-1))| \gamma \psi_M \bigg( \frac{1}{\omega}
        + t_1 \left( 1 + \rho \right) \bigg) \|\nu\|_\infty \notag
    \end{align*}
    where the second line follows from Lemmas~\ref{lem:FTpsiHSbound} and~\ref{lem:FTepsISS},
    the fourth line uses the fact that $\omega > 0$, 
    and the last line follows from the fact that $|\varepsilon(0,0)| = |\nu(0)| \leq \|\nu\|_\infty$.
    Since $t_1 \in (0, t'+T]$, we have 
    \begin{equation} \label{eqn:FTHSnuBound0122}
    \begin{aligned} 
        |\tilde \theta(t_1,1)|
        &\leq |K_1(\xi'(t_j,j-1))| \gamma \psi_M \bigg( \frac{1}{\omega} + (t'+T) \left( 1 + \rho \right) \bigg) \|\nu\|_\infty \\
        &\leq \frac{\gamma \psi_M}{1 - \sigma(\gamma,\alpha)} \bigg( \frac{1}{\omega} + (t'+T) \left( 1 + \rho \right) \bigg) \|\nu\|_\infty
        = \kappa_1'(\|\nu\|_\infty)
    \end{aligned}
    \end{equation}
    with $\kappa_1'$ as in \eqref{eqn:FTcalV}, where the last inequality follows from the definition of $K_1 = - ( \Phi - I )^{-1}$ in \eqref{eqn:FTKi1}, and we use the equivalence between the dynamics of the state component $\Phi$ in \eqref{eqn:FTHSnu2} and the state transition matrix $\Omega$ from Lemma~\ref{lem:FTPhiStability} to upper bound $|( \Phi - I )^{-1}|$ using \eqref{eqn:FTK1Bound}. 
    
    By construction, the hybrid system ${\HS}_\nu'$ jumps at most once, when $k = 0$. Since $k$ is reset to $1$ at the first jump, each solution to ${\HS}_\nu'$ flows for all hybrid time in its domain after the first jump. We have from Lemma~\ref{lem:FTHSnuSol} that, for all $t \in I^1$,
    \begin{align*} 
        \tilde \theta(t,1) 
        &= \Phi(t,1) \tilde \theta(t_1,1) - \gamma \int_{t_1}^{t} \Phi(s,1) \psi(s,1)^\top( \varepsilon(s,1) + \nu(s)) ds.
    \end{align*}
    Since $\Phi(t_1,1) = I$ and the dynamics during flows of $\Phi$ in \eqref{eqn:FTHSnu2} are equivalent to the dynamics of the state transition matrix $\Omega$ for \eqref{eqn:FTlemsysa}, it follows from \eqref{eqn:FTOmegaUpperBound} that $|\Phi(t,1)| \leq 1$ for all $t \in I^{1}$.
    Using similar arguments as in \eqref{eqn:FTHSnuBound012}, we obtain
    \begin{equation} \label{eqn:FTHSnuBound022}
    \begin{aligned}
        |\tilde \theta(t,1)| 
        %
        %
        %
        %
        %
        &\leq |\tilde \theta(t_1,1)| 
        + \gamma \psi_M \bigg( \frac{1}{\omega} \|\nu\|_\infty
        + \int_{0}^{t} \left( |\nu(s)| + \rho \|\nu\|_\infty \right) ds \bigg) \\
        &\leq \kappa_1'(\|\nu\|_\infty)
        + \kappa_2'(\|\nu\|_\infty)
        + \int_{0}^{t} \kappa_3'(|\nu(s)|,\|\nu\|_\infty) ds 
    \end{aligned}
    \end{equation}
    for all $t \in I^1$, with $\kappa_1', \kappa_2', \kappa_3'$ as in \eqref{eqn:FTcalV}, where 
    the last line follows from substituting \eqref{eqn:FTHSnuBound0122}.
    
    We unify the bounds in \eqref{eqn:FTHSnuBound012} and \eqref{eqn:FTHSnuBound022} to obtain that, for all $(t,j) \in \dom \xi'$,
    \begin{align*} 
        |\tilde \theta(t,j)|
        &\leq q(j) |\tilde \theta(0,0)|
        + (1 - q(j)) \kappa_1'(\|\nu\|_\infty) \\
        &\qquad + \kappa_2'(\|\nu\|_\infty) 
        + \int_{0}^{t} \kappa_3'(|\nu(s)|,\|\nu\|_\infty) ds
    \end{align*}
    where $q(j) \coleq 1$ if $j = 0$ and $q(j) \coleq 0$ if $j \in \nats \setminus \{0\}$. Hence, \eqref{eqn:FTxiboundiISS} holds.
    \end{proof}
}{
\begin{proof} 
Let $\xi' = (x, \hat \theta, \psi, \eta, \Phi, Q, m, \tau, k)$ be a maximal solution to ${\HS}_\nu'$ in \eqref{eqn:FTHSnu2} from $\xi'(0,0) \in {\cal X}_{\nu,0}'$. 
We have from Lemma~\ref{lem:FTHSnuSol} and from similar arguments as in \eqref{eqn:FTISSErrorBoundI0} that the parameter estimation error $\tilde \theta = \theta - \hat \theta$ satisfies
\begin{align} \label{eqn:FTHSnuBound022}
    |\tilde \theta(t,j)| 
    &\leq |\tilde \theta(t_j,j)| \\
    &\quad + \gamma \psi_M \bigg( \frac{1}{\omega} \|\nu\|_\infty
    + \int_{0}^{t} \left( |\nu(s)| + \rho \|\nu\|_\infty \right) ds \bigg) \notag \\
    &\leq \mbox{$\squeezespaces{0.5} |\tilde \theta(t_j,j)|
        + \kappa_2'(\|\nu\|_\infty)
        + \displaystyle\int_{0}^{t} \kappa_3'(|\nu(s)|,\|\nu\|_\infty) ds $} \notag
    %
\end{align}
for each $j \in \{0, 1\}$ and all $t \in I^j$, with $\kappa_2',\kappa_3'$ as in \eqref{eqn:FTcalV}. 
Since $\xi'$ is maximal and Assumptions~\ref{asm:FTphiBound} and~\ref{def:FTFE} hold, it follows from Theorem~\ref{thm:FTClosed_Loop_Stab_HS1} that $\xi'$ jumps at hybrid time $(t_1,0)$ with $t_1 \in (0, T]$.
We bound $|\tilde \theta(t_1,1)|$ by employing Lemmas~\ref{lem:FTHSnuSol},~\ref{lem:FTpsiHSbound}, and~\ref{lem:FTepsISS} and using similar arguments as in \eqref{eqn:FTHsnuJumpBound}-\eqref{eqn:FTHSnuBound01} as
\begin{align} \label{eqn:FTHSnuBound0122}
    &|\tilde \theta(t_1,1)| \\
    &\leq |K_1(\xi'(t_1,0))| \gamma \psi_M \bigg( \frac{1}{\omega} + T \left( 1 + \rho \right) \bigg) \|\nu\|_\infty \notag \\
    &\leq \mbox{$\squeezespaces{0.8} \displaystyle\frac{\gamma \psi_M}{1 - \sigma(\gamma,\alpha)} \bigg( \displaystyle\frac{1}{\omega} + T \left( 1 + \rho \right) \bigg) \|\nu\|_\infty
    = \kappa_1'(\|\nu\|_\infty) $} \notag
\end{align}
with $\kappa_1'$ as in \eqref{eqn:FTcalV}.  
By substituting \eqref{eqn:FTHSnuBound0122} into \eqref{eqn:FTHSnuBound022}, we obtain \eqref{eqn:FTxiboundiISS}.
\end{proof}
}
\begin{example}
Recall Example~\ref{exp:FTNumericalExample1}. We employ ${\HS}_\nu'$ in \eqref{eqn:FTHSnu2} to estimate the unknown parameter vector $\theta$ in \eqref{eqn:FTexamplePlant} with additive noise $t \mapsto \nu(t) = 0.5 \sin(10 t) [1 \ 1 \ 1]^\top$ in the measurements of $x$.
The simulation is performed using the same initial conditions and parameters as Example~\ref{exp:FTNumericalExample1}, and it can be shown numerically that, along the resulting solution, $\phi$ in \eqref{eqn:FTexamplePhi2} satisfies Assumption~\ref{asm:FTphiBound}, and satisfies Assumption~\ref{def:FTFE} with $t' = 0$ sec, $T = 2.97$ sec, $\mu_\phi = 1.75$, and $L_\phi = 15.6$. Moreover, $\mu$ in \eqref{eqn:FT_mubar} is greater than $0.2$, and thus $\alpha < \mu$.
We saturate the values of $x$ and $\hat \theta$ used by the control signal $\kappa$ in \eqref{eqn:FTexamplePhi} so that they remain in a compact set. Since $\kappa$ is continuously differentiable and the reference value $x_r$ is constant, the resulting control signal is bounded and globally Lipschitz in $x$, uniformly in $\hat \theta$. Hence, $f$ in \eqref{eqn:FTexamplePhi} satisfies Assumption~\ref{asm:FTboundedf}, and it follows that the conditions of Theorem~\ref{thm:FTiISS} are satisfied.
The parameter estimation error for ${\HS}_\nu'$, shown in orange in Figure~\ref{fig:FTexampleEstimationError} alongside the error from the noise-free case in blue (see Example~\ref{exp:FTNumericalExample1}), is bounded by a function of the integral of the measurement noise in accordance with Theorem~\ref{thm:FTiISS}.
\begin{figure}[!htb]
    \centering
        \includegraphics[width=1\linewidth,height=\textheight,keepaspectratio]{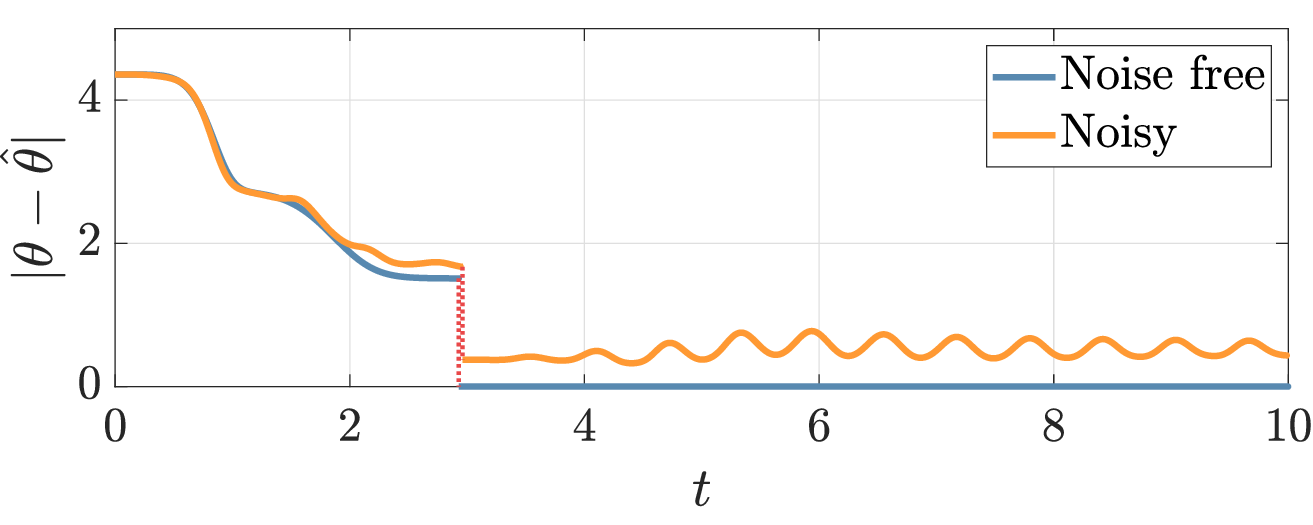}
        \caption{The projection onto $t$ of the norm of the parameter estimation error for ${\HS}_\nu'$ and $\HS$.}
    \label{fig:FTexampleEstimationError}
\end{figure}
\end{example}
\vspace{-16pt}

\section{Spacecraft Bias Torque Estimation} \label{sec:FTNumericalExample}
%

%
Consider the problem of estimating a constant disturbance torque applied to a spacecraft that is controlled by reaction wheels (RW). 
Such bias torques may arise in practice due to aerodynamic effects, gravity gradients, or solar radiation pressure differentials.
We model the spacecraft as a rigid body with state $x = (q, \omega)$, where $q \in \reals^4$ is a quaternion representing the spacecraft pointing angle, and $\omega \in \reals^3$ is the spacecraft angular velocity.
For simplicity, we consider only the rotational dynamics of the spacecraft, given by \cite{Wie1998SpacecraftDynamics}
\begin{equation} \label{eqn:PESC1}
    \begin{aligned} 
        \dot q &= \frac{1}{2} \omega_q \otimes q  \\
        \dot \omega &= J^{-1} ( - \omega \times J \omega + \kappa + \theta )
    \end{aligned}
\end{equation}
where $\omega_q \coleq (\omega, 0)$, $J \in \reals^3$ is the known positive definite spacecraft inertia matrix, $\kappa \in \reals^3$ is the known RW torque, and $\theta \in \reals^3$ is an unknown constant bias torque.
%
%

When the bias torque is nonzero, the industry-standard proportional-derivative (PD) control scheme for the RW torque fails to yield zero pointing error in steady-state. In this case, a feedforward term is added that compensates for the effect of the bias torque using an estimate of the bias, denoted by $\hat \theta$. Given $u \coleq (q_{\rm des}, \omega_{\rm des}, \hat \theta)$, where $q_{\rm des} \in \reals^4$ is a desired spacecraft pointing angle and $\omega_{\rm des} \in \reals^3$ is a desired angular velocity, the RW torque is \cite{Wie1998SpacecraftDynamics}
\begin{equation} \label{eqn:PEu}
    \kappa(x, u) \coleq - K_P \, q_{v,\rm err} - K_D \, \omega_{\rm err} - \hat \theta,
\end{equation}
where the positive definite matrices $K_P, K_D \in \reals^{3 \times 3}$ are design parameters, $\omega_{\rm err} \coleq \omega - \omega_{\rm des}$, and $q_{v,\rm err} \coleq (q_{1, \rm err}, q_{2, \rm err}, q_{3, \rm err})$ with 
\begin{equation*}
(q_{1, \rm err}, q_{2, \rm err}, q_{3, \rm err}, q_{4, \rm err}) = q \otimes q_{\rm des}^{-1}.
\end{equation*}
\begin{remark}
The desired spacecraft pointing angle can be maintained only if an equivalent RW torque is delivered to counteract the bias torque. If the bias torque is nonzero, the angular velocity of the RW constantly increases to counteract the disturbance and the RW motor eventually reaches its maximum angular velocity. 
In order to avoid the RW motor from becoming saturated, ``momentum dumping'' is applied to decrease the angular velocity of the RW \cite{Wie1998SpacecraftDynamics}. This procedure involves firing reaction control thrusters to generate a torque that is compensated by the attitude controller by actions that cause the RW to reduce their angular momentum. We neglect analyzing the effect of the momentum dumping in this example, which is reasonable since the proposed estimator achieves finite-time convergence of the parameter estimate and can therefore be employed between such events.
\end{remark}

\begin{figure}[!htb]
    \centering
        \includegraphics[width=1\linewidth,height=\textheight,keepaspectratio]{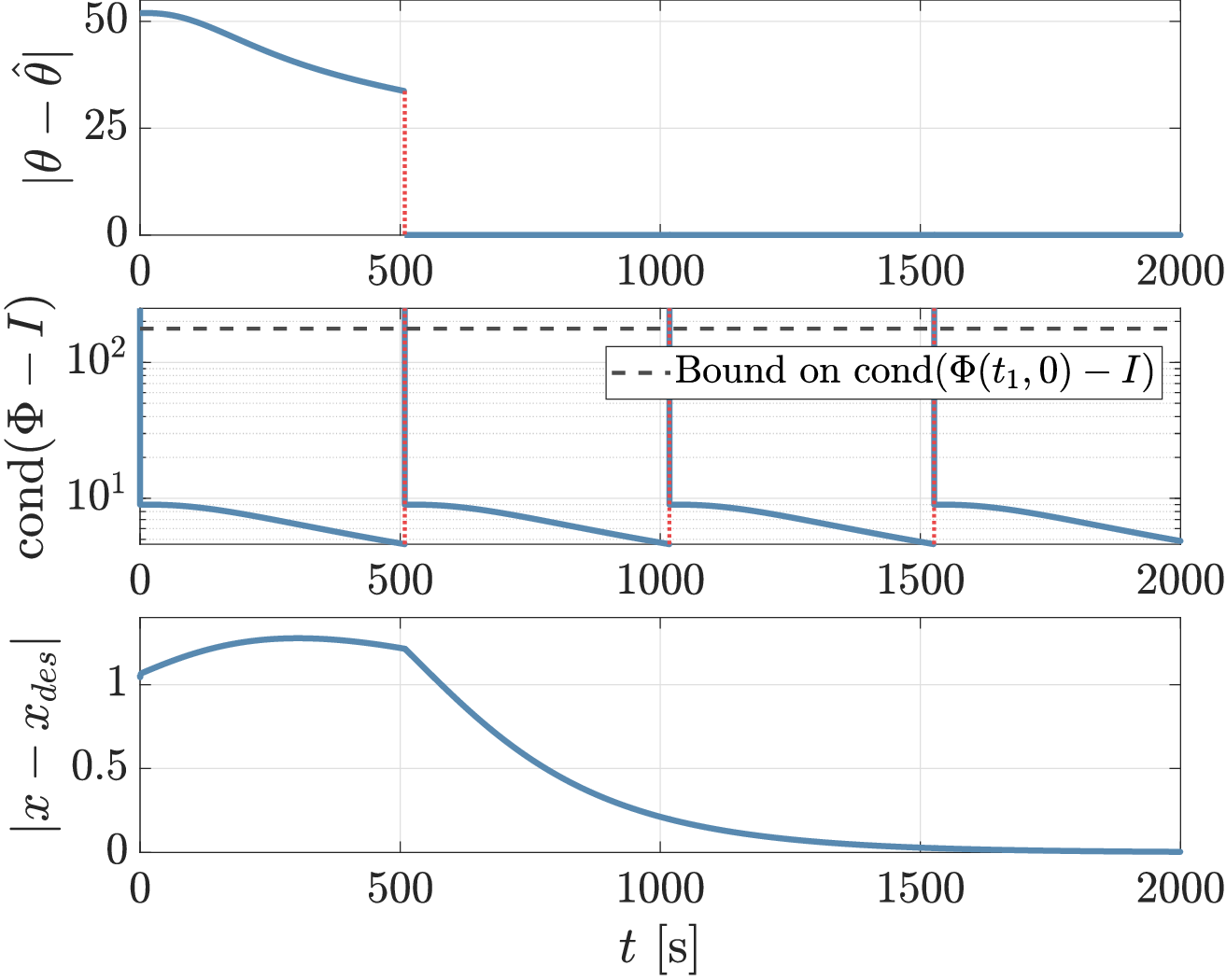}
        \caption{Top: The projection onto $t$ of the norm of the bias torque estimation error for $\HS$. 
        Middle: The projection onto $t$ of cond$(\Phi - I)$ in blue along with the upper bound on cond$(\Phi(t_1,0) - I)$ from \eqref{eqn:FTdeltabound} in dashed black.
        Bottom: The norm of the spacecraft pointing and angular velocity error.}
    \label{fig:FTSpacecraft}
\end{figure}

%
%
To estimate the unknown bias torque using the proposed algorithm, we express the spacecraft dynamics \eqref{eqn:PESC1}-\eqref{eqn:PEu} in the form of \eqref{eqn:FTxdot} with
\begin{equation} \label{eqn:FTSCmodel}
\!\!\!
\mbox{$\squeezespaces{0.6}
\allowdisplaybreaks
\begin{aligned}
    f(x,u) \coleq
    \matt{\displaystyle\frac{1}{2} \, \omega_q \otimes q \\
          J^{-1} ( - \omega \times J \omega + \kappa(x, u) )},
    \phi \coleq
    \matt{0 \\ J^{-1}}
\end{aligned}
$}
\end{equation}
%
%
%
%
We simulate $\HS$ in \eqref{eqn:FTH1} with $f$ and $\phi$ in \eqref{eqn:FTSCmodel}. 
The spacecraft system has parameters
$q_{\rm des} = (-0.5, -0.5, -0.5, 0.5)$, $\omega_{\rm des} = (0, 0, 0)$ rad/sec, $J =  \operatorname{diag}(3000, 2500, 1000)$ kg-m$^2$, $K_p = 40 \, I$, $K_d = 5000 \, I$, and unknown bias torque $\theta = (30, 30, 30)$ N-m. The initial condition for $q$ is $(0, 0, 0, 1)$ and $\omega$ is $(0.2, 0.2, 0.2)$ rad/sec.
The estimator has design parameters $\gamma = 0.5$, $\lambda = 0.01$, and $\alpha = 0.4$, and initial conditions $\hat \theta = m = 0$, $\psi = 0$, $\eta = -x$, $\Phi = I$, and $Q = 0$.
With these parameters and initial conditions, it can be shown\footnote{Code at \href{https://github.com/HybridSystemsLab/HybridFT_SpacecraftBiasTorque}{\url{https://github.com/HybridSystemsLab/HybridFT_SpacecraftBiasTorque}}} that $\phi$ in \eqref{eqn:FTSCmodel} satisfies Assumption~\ref{asm:FTphiBound} with $\phi_M = 1.13 \times 10^{-3}$, and satisfies Assumption~\ref{def:FTFE} with $t' = 0$ sec, $T = 509$ sec, $\mu_\phi = 5.65 \times 10^{-5}$, and $L_\phi = 0$. Moreover, $\mu$ in \eqref{eqn:FT_mubar} is greater than $0.4$, and thus $\alpha < \mu$. Hence, the conditions of Theorem~\ref{thm:FTClosed_Loop_Stab_HS1} are satisfied.
The bias torque estimation error converges to zero in finite time in accordance with Theorem \ref{thm:FTClosed_Loop_Stab_HS1}, as shown in the top sub-plot of Figure \ref{fig:FTSpacecraft}.
The middle subplot shows the condition number of $\Phi - I$ in blue along with the bound \eqref{eqn:FTdeltabound} in dashed black. In this case, \eqref{eqn:FTdeltabound} ensures that cond$(\Phi(t_1,0) - I) \leq 177.3$ and the actual value is cond$(\Phi(t_1,0) - I) = 4.6$.
Following convergence of the bias torque estimate, the spacecraft state $x = (q, \omega)$ converges to the desired value $x_{\rm des} = (q_{\rm des}, \omega_{\rm des})$, as shown in the bottom sub-plot of Figure~\ref{fig:FTSpacecraft}.
\vspace{-16pt}
%

\section{Conclusion} \label{sec:FTConclusion}
In this paper, we proposed a hybrid algorithm for finite-time estimation of unknown constant parameters in a class of continuous-time systems with nonlinear dynamics that are affine in the unknown parameters. We showed that the algorithm guarantees finite-time convergence of the parameter estimate to the true value under excitation over only a finite interval.
Furthermore, we showed that the parameter estimate is numerically well conditioned at the convergence time by upper bounding the condition number of a matrix whose inverse is used to reset the parameter estimate at jumps. 
We showed that, under persistence of excitation, the parameter estimate is input-to-state stable with respect to noise in the measurements of the plant state, and that, under excitation on a finite interval, the parameter estimate is bounded by a function of the integral of the measurement noise.
Future work on this topic includes extending the proposed algorithm to estimate time-varying unknown parameters, and developing finite-time estimation algorithms for systems whose dynamics depend nonlinearly on the unknown parameters.
%
\vspace{-8pt}

\medskip
\bibliographystyle{unsrt}
\bibliography{Bibliography}

\begin{thebibliography}{10}

\bibitem{8855027}
M.~S. Rafaq and J.-W. Jung.
\newblock A comprehensive review of state-of-the-art parameter estimation
  techniques for permanent magnet synchronous motors in wide speed range.
\newblock {\em IEEE Transactions on Industrial Informatics}, 16(7):4747--4758,
  2020.

\bibitem{8466598}
K.~Dehghanpour, Z.~Wang, J.~Wang, Y.~Yuan, and F.~Bu.
\newblock A survey on state estimation techniques and challenges in smart
  distribution systems.
\newblock {\em IEEE Transactions on Smart Grid}, 10(2):2312--2322, 2019.

\bibitem{doi:10.2514/3.20454}
K.~W. Iliff.
\newblock Parameter estimation for flight vehicles.
\newblock {\em Journal of Guidance, Control, and Dynamics}, 12(5):609--622,
  1989.

\bibitem{Narendra_adaptive_1989}
K.~S. Narendra and A.~M. Annaswamy.
\newblock {\em {Stable Adaptive Systems}}.
\newblock Prentice-Hall, Inc., USA, 1989.

\bibitem{Sastry1989Adaptive}
S.~Sastry and M.~Bodson.
\newblock {\em {Adaptive Control: Stability, Convergence, and Robustness}}.
\newblock Prentice-Hall, Inc., USA, 1989.

\bibitem{Khalil2002Nonlinear}
H.~K Khalil.
\newblock {\em {Nonlinear Systems; 3rd ed.}}
\newblock Prentice-Hall, Upper Saddle River, NJ, 2002.

\bibitem{243}
A.~Saoud and R.~G. Sanfelice.
\newblock {A Robust Hybrid Finite-Time Parameter Estimator With Relaxed
  Persistence of Excitation Condition}.
\newblock In {\em 2021 American Control Conference (ACC)}, pages 3975--3980,
  2021.

\bibitem{65}
R.~Goebel, R.~G. Sanfelice, and A.~R. Teel.
\newblock {\em {Hybrid Dynamical Systems: Modeling, Stability, and
  Robustness}}.
\newblock Princeton University Press, New Jersey, 2012.

\bibitem{220}
R.~G. Sanfelice.
\newblock {\em {Hybrid Feedback Control}}.
\newblock Princeton University Press, New Jersey, 2021.

\bibitem{Wang2019FiniteTime}
J.~Wang, D.~Efimov, and A.A. Bobtsov.
\newblock Finite-time parameter estimation without persistence of excitation.
\newblock In {\em 2019 18th European Control Conference (ECC)}, pages
  2963--2968, 2019.

\bibitem{Ortega1988LeastSquares}
R.~Ortega.
\newblock An on-line least-squares parameter estimator with finite convergence
  time.
\newblock {\em Proceedings of the IEEE}, 76(7):847--848, 1988.

\bibitem{Adetola2007FiniteTime}
V.~Adetola and M.~Guay.
\newblock Finite-time parameter estimation in adaptive control of nonlinear
  systems.
\newblock In {\em 2007 American Control Conference}, pages 1920--1924, 2007.

\bibitem{Na2015FiniteTime}
J.~Na, M.~N. Mahyuddin, G.~Herrmann, X.~Ren, and P.~Barber.
\newblock Robust adaptive finite-time parameter estimation and control for
  robotic systems.
\newblock {\em International Journal of Robust and Nonlinear Control},
  25(16):3045--3071, 2015.

\bibitem{Vahidi2021Memory}
A.~Vahidi-Moghaddam, M.~Mazouchi, and H.~Modares.
\newblock Memory-augmented system identification with finite-time convergence.
\newblock {\em IEEE Control Systems Letters}, 5(2):571--576, 2021.

\bibitem{Mohamed2026DREM}
M.~Kharrat.
\newblock {DREM}-based predefined-time adaptive learning control for nonlinear
  systems with actuator faults and input delay.
\newblock {\em European Journal of Control}, 87:101426, 2026.

\bibitem{Ortega2022ParamId}
R.~Ortega, A.~Bobtsov, and N.~Nikolaev.
\newblock Parameter identification with finite-convergence time alertness
  preservation.
\newblock {\em IEEE Control Systems Letters}, 6:205--210, 2022.

\bibitem{Saoud2024HybridPersistency}
Adnane Saoud, Mohamed Maghenem, Antonio Loría, and Ricardo~G. Sanfelice.
\newblock Hybrid persistency of excitation in adaptive estimation for hybrid
  systems.
\newblock {\em IEEE Transactions on Automatic Control}, 69(12):8828--8835,
  2024.

\bibitem{242}
A.~Saoud, M.~Maghenem, and R.~G. Sanfelice.
\newblock {A Hybrid Gradient Algorithm for Linear Regression with Hybrid
  Signals}.
\newblock In {\em 2021 American Control Conference (ACC)}, pages 4997--5002,
  2021.

\bibitem{256}
R.~S. Johnson, A.~Saoud, and R.~G. Sanfelice.
\newblock {Robust Finite-Time Parameter Estimation for Linear Dynamical
  Systems}.
\newblock In {\em 2021 60th IEEE Conference on Decision and Control (CDC)},
  pages 4654--4659, 2021.

\bibitem{SaoudFiniteTimeInPrep}
A.~Saoud, R.~S. Johnson, and R.~G. Sanfelice.
\newblock Robust {Hybrid} {Finite} {Time} {Parameter} {Estimation} {Without}
  {Persistence} {of} {Excitation}.
\newblock {\em IEEE Transactions on Automatic Control}, Under review.

\bibitem{hartman2012robust}
M.~{Hartman}, N.~{Bauer}, and A.~R. {Teel}.
\newblock {Robust Finite-Time Parameter Estimation Using a Hybrid Systems
  Framework}.
\newblock {\em IEEE Transactions on Automatic Control}, 57(11):2956--2962,
  2012.

\bibitem{Panteley2002Detectability}
E.~Panteley, R.~Ortega, and P.~Moya.
\newblock {Overcoming the detectability obstacle in certainty equivalence
  adaptive control}.
\newblock {\em Automatica}, 38(7):1125--1132, 2002.

\bibitem{Adetola_finite_time_est_2008}
V.~{Adetola} and M.~{Guay}.
\newblock {Finite-Time Parameter Estimation in Adaptive Control of Nonlinear
  Systems}.
\newblock {\em IEEE Transactions on Automatic Control}, 53(3):807--811, 2008.

\bibitem{Jung-Shan1999Nonlinearities}
J.-S. Lin and I.~Kanellakopoulos.
\newblock {Nonlinearities enhance parameter convergence in strict feedback
  systems}.
\newblock {\em IEEE Transactions on Automatic Control}, 44(1):89--94, 1999.

\bibitem{Horn2012Matrix}
R.~A. Horn and C.~R. Johnson.
\newblock {\em {Matrix Analysis}}.
\newblock Cambridge University Press, USA, 2nd edition, 2012.

\bibitem{Wie1998SpacecraftDynamics}
B.~Wie.
\newblock {\em Space vehicle dynamics and control}.
\newblock AIAA education series. American Institute of Aeronautics and
  Astronautics, Reston, VA, 1998.

\end{thebibliography}

\appendix
\iftoggle{long}{\appendix}{\begin{appendices}}
\renewcommand{\thehelptheorem}{\Alph{section}.\arabic{helptheorem}}

%
\vspace{-10pt}
\section{Proof of Lemma~\ref{lem:FT_FEPsi}} \label{apx:FT_FEPsi}
\vspace{-6pt}
\begin{proof}[\unskip\nopunct]
We follow an approach inspired by the proof of \cite[Lemma~2.6.7]{Sastry1989Adaptive}.
Let $\psi : [t', \infty) \to \reals^{n \times p}$ be the solution to \eqref{eqn:FTpsidynamics0} and
pick $c \in \reals^p$ with $|c| = 1$. For each $s \in [t', t'+T]$, define
\vspace{-6pt}
\begin{align} \label{eqn:FT_uyDef}
    y(s) &\coleq \psi(s) \, \! c, &
    u(s) &\coleq \phi(s) \, \! c.
\end{align}
From Assumption~\ref{asm:FTphiBound} and from \eqref{eqn:Lischitz_phi} in Assumption~\ref{def:FTFE}, we have that, for all $s, v \in [t', t' + T]$,
\begin{equation} \label{eqn:FT_uProp}
\mbox{$\squeezespaces{0.7}
\begin{aligned} 
    \!\!\!\!\!
    |u(s)| \leq |\phi(s)|_F \leq \phi_M, \quad
    |u(s) - u(v)| \leq L_\phi |s - v|
\end{aligned}
$}
\end{equation}
Moreover, from \eqref{eqn:FTFE} in Assumption~\ref{def:FTFE},
\begin{equation} \label{eqn:FT_muPhiLb}
\mbox{$\squeezespaces{0.4}
\begin{aligned} 
    \!\!\!\!\!\!\!\!
    \int_{t'}^{t'+T} |u(s)|^2 ds = c^\top \bigg( \int_{t'}^{t'+T} \phi(s)^\top \phi(s) ds \bigg) c \geq \mu_\phi.
\end{aligned}
$}
\end{equation}
To relate the excitation of $u$ to that of $y$ in \eqref{eqn:FT_uyDef}, we define 
\begin{equation} \label{eqn:FT_rDef}
    r(s) \coleq u(s) - ay(s)
\end{equation}
for all $s \in [t', t'+T]$, with $a$ as in \eqref{eqn:FT_mu_a}. 
By \eqref{eqn:FT_uProp}, $u$ in \eqref{eqn:FT_uyDef} is Lipschitz on $[t', t'+T]$, and thus differentiable for almost all $s \in [t', t'+T]$. Hence, using that 
\iftoggle{report}{}{\linebreak}
$\dot y = -\lambda y + u(s)$, for almost all $s \in [t', t'+T]$,
\begin{equation} \label{eqn:FT_rdyn}
\mbox{$\squeezespaces{0.5}
\begin{aligned}
    \!\!\!\!\!\!\!
    \dot r = \dot u - a (-\lambda y + u) = -a r + h, \quad
    h \coleq \dot u - a (a - \lambda) y
\end{aligned}
$}
\end{equation}
where the argument $s$ of $r, u, y, h$ is omitted for readability. 
Consider the Lyapunov function
\begin{align} \label{eqn:FT_Wdef}
    W(r) &\coleq \frac{1}{2} r^\top r \qquad \forall r \in \reals^n.
\end{align}
Differentiating $W$ along trajectories of \eqref{eqn:FT_rdyn} and using that, for any $v, w \in \reals^n$, $v^\top w \leq \frac{a}{2}|v|^2 + \frac{1}{2a}|w|^2$, we obtain
\begin{align} \label{eqn:FT_dotW}
    \frac{d}{ds} W(r(s))
    &= -a |r(s)|^2 + r(s)^\top h(s) \notag  \\
    &\leq -\frac{a}{2} |r(s)|^2 + \frac{1}{2a} |h(s)|^2
\end{align}
for almost all $s \in [t', t'+T]$. 
Since, by Assumption~\ref{def:FTFE}, $\phi$ is continuous on $[t', t'+T]$, the solution $\psi$ to \eqref{eqn:FTpsidynamics0}, and thus $y$ in \eqref{eqn:FT_uyDef}, is continuously differentiable on $[t', t'+T]$.
Moreover, by \eqref{eqn:FT_uProp}, $u$ in \eqref{eqn:FT_uyDef} is absolutely continuous on $[t', t'+T]$. 
Hence, $r$ in \eqref{eqn:FT_rDef} is absolutely continuous on $[t', t'+T]$, and so is $s \mapsto W(r(s))$ in \eqref{eqn:FT_Wdef}.
Integrating \eqref{eqn:FT_dotW} between times $t'$ and $t' + T$, we obtain
\begin{align*}
    &W(r(t'+T)) - W(r(t')) \\
    &\qquad \leq -\frac{a}{2} \int_{t'}^{t'+T} |r(s)|^2 ds
    + \frac{1}{2 a} \int_{t'}^{t'+T} |h(s)|^2 ds.
\end{align*}
Since $\psi(t') = 0$ by \eqref{eqn:FTpsidynamics0}, it follows from \eqref{eqn:FT_uyDef} and \eqref{eqn:FT_rDef} that $y(t') = 0$ and $r(t') = u(t')$. Hence, by \eqref{eqn:FT_uProp}, $W(r(t')) = \frac{1}{2}|u(t')|^2 \leq \frac{1}{2}\phi_M^2$.
Rearranging terms and using that $W(r(t'+T)) \geq 0$ by \eqref{eqn:FT_Wdef}, we obtain
\begin{align} \label{eqn:FT_intrUb}
    \int_{t'}^{t'+T} |r(s)|^2 ds \leq \frac{\phi_M^2}{a} + \frac{1}{a^2} \int_{t'}^{t'+T} |h(s)|^2 ds.
\end{align}
We bound the second term on the right-hand side of \eqref{eqn:FT_intrUb} using \eqref{eqn:FT_rdyn}, Minkowski's inequality, and the fact that $a \geq \lambda$ and, from \eqref{eqn:FT_uProp}, $|\dot u(s)| \leq L_\phi$ for almost all $s \in [t', t'+T]$, to obtain
\begin{align} \label{eqn:FT_intfUb}
    &\bigg( \int_{t'}^{t'+T} |h(s)|^2 ds \bigg)^{1/2} \notag \\
    &\leq \mbox{$\squeezespaces{0.3}
    \bigg( \displaystyle\int_{t'}^{t'+T} |\dot u(s)|^2 ds \bigg)^{1/2} 
    + a( a-\lambda ) \bigg( \displaystyle\int_{t'}^{t'+T} |y(s)|^2 ds \bigg)^{1/2} 
    $} \notag \\
    &\leq L_\phi \sqrt{T} 
    + a( a-\lambda ) \bigg( \int_{t'}^{t'+T} |y(s)|^2 ds \bigg)^{1/2}.
\end{align}
Taking the square root of both sides of \eqref{eqn:FT_intrUb}, substituting \eqref{eqn:FT_intfUb}, and using that $\sqrt{b_1 + b_2} \leq \sqrt{b_1} + \sqrt{b_2}$ for all $b_1, b_2 \geq 0$, we have
\begin{align} \label{eqn:FT_rUb}
    &\bigg( \int_{t'}^{t'+T} |r(s)|^2 ds \bigg)^{1/2} \\
    &\quad \leq \frac{\phi_M}{\sqrt{a}} + \frac{L_\phi \sqrt{T}}{a} 
    + ( a-\lambda ) \bigg( \int_{t'}^{t'+T} |y(s)|^2 ds \bigg)^{1/2}. \notag
\end{align}
We now return to \eqref{eqn:FT_muPhiLb}. Since, by \eqref{eqn:FT_rDef}, $u = ay + r$, it follows from Minkowski's inequality that
\begin{align} \label{eqn:FT_uUb}
    \sqrt{\mu_\phi} 
    &\leq \bigg( \int_{t'}^{t'+T} |u(s)|^2 ds \bigg)^{1/2} \\
    &\leq a \bigg( \int_{t'}^{t'+T} |y(s)|^2 ds \bigg)^{1/2} 
    + \bigg( \int_{t'}^{t'+T} |r(s)|^2 ds \bigg)^{1/2}. \notag
\end{align}
Substituting \eqref{eqn:FT_rUb} into \eqref{eqn:FT_uUb}, we obtain
\begin{align*}
    \sqrt{\mu_\phi}
    \leq (2a - \lambda) \bigg( \int_{t'}^{t'+T} |y(s)|^2 ds \bigg)^{1/2}
    + \frac{\phi_M}{\sqrt{a}} + \frac{L_\phi \sqrt{T}}{a}.
\end{align*}
From the definition of $a$ in \eqref{eqn:FT_mu_a}, 
$\phi_M / \sqrt{a} \leq \sqrt{\mu_\phi} / 4$ 
and
$L_\phi \sqrt{T} / a \leq \sqrt{\mu_\phi} / 4$.
Hence,
\begin{align*}
    \frac{\sqrt{\mu_\phi}}{2}
    \leq (2a - \lambda) \bigg( \int_{t'}^{t'+T} |y(s)|^2 ds \bigg)^{1/2}
\end{align*}
and, since $2a - \lambda > 0$,
\begin{align*}
    \bigg( \int_{t'}^{t'+T} |y(s)|^2 ds \bigg)^{1/2}
    \geq \frac{\sqrt{\mu_\phi}}{2(2a - \lambda)}.
\end{align*}
Squaring both sides and using the definition of $y$ in \eqref{eqn:FT_uyDef}, we obtain
\begin{align*}
    c^\top \bigg( \int_{t'}^{t'+T} \psi(s)^\top \psi(s) ds \bigg) c
    &= \int_{t'}^{t'+T} |y(s)|^2 ds \\
    &\geq \frac{\mu_\phi}{4 (2a - \lambda)^2}
    = \mu_\psi.
\end{align*}
Since this inequality holds for any $c \in \reals^p$ with $|c| = 1$, we conclude that \eqref{eqn:FT_psiFE} holds with $\mu_\psi$ as in \eqref{eqn:FT_mu_a}.
\end{proof}

\section{Proof of Lemma~\ref{lem:FTHSnuSol}} \label{apx:FTHSnuSol}
\begin{proof}[\unskip\nopunct]
Let $\xi$ be a solution to $\HS_\nu$ from $\xi(0,0) \in {\cal X}_{\nu,0}$, where ${\cal X}_{\nu,0}$ is defined above \eqref{eqn:FTxiboundISS}.
We first analyze the evolution of the parameter estimation error $(t,j) \mapsto \tilde \theta(t,j) \coleq \theta - \hat \theta(t,j)$ during flows. Since $\theta$ is constant, it follows from $F_\nu$ in \eqref{eqn:FTFnuGnu} that
\begin{align*}
    &\frac{d}{dt}{\tilde \theta}(t,j)
     = - \frac{d}{dt}{\hat \theta}(t,j) \\
    &= - \gamma \psi(t,j)^\top \big( x(t,j) + \nu(t) + \eta(t,j) - \psi(t,j) \hat \theta(t,j) \big) \\
    &= - \gamma \psi(t,j)^\top \psi(t,j) \tilde \theta(t,j) - \gamma \psi(t,j)^\top \big( \varepsilon(t,j) + \nu(t) \big)
\end{align*}
%
%
%
for each $j \in \nats$ and all $t \in I^j$, with $\varepsilon$ as in \eqref{eqn:FTeps}. 
Since $\Phi(0,0) = I$ and, by the jump map $G_\nu$ in \eqref{eqn:FTFnuGnu}, $\Phi$ is equal to the identity immediately after each jump, it follows from $F_\nu$ in \eqref{eqn:FTFnuGnu} that, since $\Phi$ is the state transition matrix for \eqref{eqn:FTtildethetaC},
\begin{align*}
    \tilde \theta(t,j) 
    &= \Phi(t,j) \tilde \theta(t_j,j) \\ 
    &\qquad - \gamma \int_{t_j}^{t} \Phi(s,j) \psi(s,j)^\top( \varepsilon(s,j) + \nu(s)) ds
\end{align*}
for each $j \in \nats$ and all $t \in I^j$.
Thus, \eqref{eqn:FTTildeThetaFlow} holds.

Next, we analyze the evolution of $\tilde \theta$ at jumps. 
Let $(t,j) \in \dom \xi$ be such that $(t,j+1) \in \dom \xi$. 
Denoting $\tilde m = \theta - m$ and omitting the argument $\xi$ of $K_1$ for readability, it follows from $G_\nu$ in \eqref{eqn:FTFnuGnu} that
\begin{equation*}
    \begin{aligned}
        &\tilde \theta(t,j+1) \\
        &= \theta - R(\xi(t,j)) 
        = \tilde m(t,j) + K_1 \big( \tilde \theta(t,j) - \tilde m(t,j) \big) \\
        &= \tilde m(t_j,j) + K_1 \big( \Phi(t,j) - I \big) \tilde m(t_j,j) \\
        &\quad - K_1 \gamma \int_{t_j}^{t} \Phi(s,j) \psi(s,j)^\top \big( \varepsilon(s,j) + \nu(s) \big) ds \\
        &= - K_1 \gamma \int_{t_j}^{t} \Phi(s,j) \psi(s,j)^\top \big( \varepsilon(s,j) + \nu(s) \big) ds
    \end{aligned}
\end{equation*} 
where the first line follows from the definition of $R$ in \eqref{eqn:FTR}, 
the second line follows from substituting \eqref{eqn:FTTildeThetaFlow} 
and from the fact that, since $\xi(0,0) \in {\cal X}_{\nu,0}$, $\tilde m = \tilde \theta$ at the beginning of each interval of flow in $\xi$ and $m$ is constant during flows, 
and the third line follows from the definition of $\K_1$ in \eqref{eqn:FTKi1}.
Thus, \eqref{eqn:FTTildeThetaJump} holds.
\end{proof}
\section{Proof of Lemma~\ref{lem:FTepsISS}} \label{apx:FTepsISS}
\begin{proof}[\unskip\nopunct]
Consider the Lyapunov function
\begin{align} \label{eqn:FTVepsilon}
    V_\varepsilon(\xi) 
    &\coleq \frac{1}{2} \varepsilon^\top \varepsilon
    \quad \forall \xi \in C_\nu \cup D_\nu
\end{align}
with $\varepsilon$ as in \eqref{eqn:FTeps}. Since $\theta$ is constant, we have from the definition of $F_{\nu}$ below \eqref{eqn:FTHSnu} that, during flows,
\begin{align*}
    \dot \varepsilon 
    &= \dot x + \dot \eta - \dot \psi \theta \\
    &= - \lambda \varepsilon + f(x,u(\tau)) - f(x+\nu(\tau),u(\tau)) - \lambda \nu(\tau).
\end{align*}
Hence, for all $\xi \in C_\nu$,
\begin{align*} %
    &\langle \nabla V_\varepsilon(\xi), F_\nu(\xi) \rangle \\
    &= - \lambda \varepsilon^\top \varepsilon + \varepsilon^\top( f(x,u(\tau))
    - f(x+\nu(\tau),u(\tau)) - \lambda \nu(\tau) ) \\
    &\leq - \lambda \varepsilon^\top \varepsilon + |\varepsilon| (L_x + \lambda) |\nu(\tau)|
    \leq - 2 \lambda V_\varepsilon(\xi) + |\varepsilon| \chi(\tau)
\end{align*}
where the first inequality follows from Assumption~\ref{asm:FTboundedf}, and we define $\chi(\tau) \coleq (L_x + \lambda) |\nu(\tau)|$ for readability. 
We use the fact that, via Young's inequality, for any $\varrho > 0$, $|\varepsilon| \chi(\tau) \leq \varrho |\varepsilon|^2 + \varrho^{-1} \chi(\tau)^2$ for all $\tau \geq 0$. Choosing $\varrho = \lambda / 2$,
\begin{align*}
    \langle \nabla V_\varepsilon(\xi), F_\nu(\xi) \rangle 
    &\leq - \lambda V_\varepsilon(\xi) + \frac{2}{\lambda} \chi(\tau)^2.
\end{align*}
Let $\zeta \in (0, 1)$. 
We rewrite the expression above as
\begin{equation*}
    \mbox{$\squeezespaces{0.9}
    \langle \nabla V_\varepsilon(\xi), F_\nu(\xi) \rangle 
    \leq - \lambda (1 - \zeta) V_\varepsilon(\xi) - \lambda \zeta V_\varepsilon(\xi) + \displaystyle\frac{2}{\lambda} \chi(\tau)^2.
    $}
\end{equation*}
Hence,
\begin{align} \label{eqn:FTetaISSflow}
    &\langle \nabla V_\varepsilon(\xi), F_\nu(\xi) \rangle
    \leq - \lambda (1 - \zeta) V_\varepsilon(\xi) \\
    &\qquad \qquad \forall \xi \in C_\nu : \ V_\varepsilon(\xi) \geq \frac{2 (L_x + \lambda)^2}{\lambda^2 \zeta} |\nu(\tau)|^2. \notag
\end{align}
Let us now analyze the variation of $V_\varepsilon$ at jumps. Since $\theta$ is constant, we have from the definition of $G_{\nu}$ below \eqref{eqn:FTHSnu} that, at jumps,
%
$
    \varepsilon^+ 
    = x^+ + \eta^+ - \psi^+ \theta
    = - \nu(\tau)
$.
%
Hence, for all $\xi \in D_\nu$,
\begin{equation*}
\begin{aligned} 
    V_\varepsilon(G_\nu(\xi)) - V_\varepsilon(\xi)
    = - V_\varepsilon(\xi) + \frac{1}{2} |\nu(\tau)|^2
\end{aligned}
\end{equation*}
Let $\zeta \in (0, 1)$. 
We rewrite the expression above as
\begin{equation*}
    \mbox{$\squeezespaces{0.9}
    V_\varepsilon(G_\nu(\xi)) - V_\varepsilon(\xi)
    = - (1 - \zeta) V_\varepsilon(\xi) - \zeta V_\varepsilon(\xi) + \displaystyle\frac{1}{2} |\nu(\tau)|^2
    $}
\end{equation*}
Hence,
\begin{align} \label{eqn:FTetaISSjump}
    &V_\varepsilon(G_\nu(\xi)) - V_\varepsilon(\xi)
    \leq - (1 - \zeta) V_\varepsilon(\xi), \\
    &\qquad \qquad \qquad \qquad \quad
    \forall \xi \in D_\nu : \ V_\varepsilon(\xi) \geq \frac{1}{2 \zeta} |\nu(\tau)|^2. \notag
\end{align}
Combining the expressions in \eqref{eqn:FTetaISSflow} and \eqref{eqn:FTetaISSjump}, we obtain
\begin{align*}
    \langle \nabla V_\varepsilon(\xi), F_\nu(\xi) \rangle 
    &\leq - \lambda (1 - \zeta) V_\varepsilon(\xi) &
    \forall \xi &\in C_\nu \cap S_\varepsilon \\
    V_\varepsilon(G_\nu(\xi)) - V_\varepsilon(\xi)
    &\leq - (1 - \zeta) V_\varepsilon(\xi) &
    \forall \xi &\in D_\nu \cap S_\varepsilon
\end{align*}
where
\begin{align*}
    \mbox{$\squeezespaces{0.5}
    S_\varepsilon \coleq \left\{\xi \in C_\nu \cup D_\nu : V_\varepsilon(\xi) \geq \max \left\{
        \frac{2 (L_x + \lambda)^2}{\lambda^2 \zeta}, 
        \frac{1}{2 \zeta}
        \right\}
        |\nu(\tau)|^2
    \right\}.
    $}
\end{align*}
\noindent Then, for each solution $\xi$ to $\HS_\nu$, by integration using the bounds above, we obtain
\begin{align*}
    &V_\varepsilon(\xi(t,j)) \\
    &\leq \exp \big\{ 
    -\lambda (1 - \zeta) (t - t')
    + \ln \left( \zeta \right) (j - j')    
    \big\} V_\varepsilon(\xi(t',j')) \\
    &\qquad+ 
    \max \left\{
        \frac{2 (L_x + \lambda)^2}{\lambda^2 \zeta}, 
        \frac{1}{2 \zeta}
        \right\}
        \|\nu\|_\infty^2
\end{align*}
for all $(t,j), (t', j') \in \dom \xi$ satisfying $t+j \geq t'+j'$.
Using the definition of $V_\varepsilon$ in \eqref{eqn:FTVepsilon}, we conclude that \eqref{eqn:FTepsISS} holds with $\omega$ in \eqref{eqn:FTomega}.
\end{proof}

\iftoggle{long}{}{\end{appendices}}

\end{document}